\documentclass[a4paper]{amsart}
\usepackage{amsmath,amssymb}
\usepackage{graphicx}
\usepackage[all]{xy}

\usepackage{multirow}
\newtheorem{Theorem}{Theorem}[section]
\newtheorem{Lemma}[Theorem]{Lemma}
\newtheorem{Proposition}[Theorem]{Proposition}
\newtheorem{Corollary}[Theorem]{Corollary}

\newtheorem{Example}[Theorem]{Example}

\newtheorem{Remark}[Theorem]{Remark}

\makeatletter
\@addtoreset{figure}{section}
\def\@thmcountersep{-}
\makeatother

\numberwithin{equation}{section}

\begin{document}

\title[Crossing numbers and rotation numbers of cycles in a graph]{Crossing numbers and rotation numbers of cycles in a plane immersed graph}

\author{Ayumu Inoue}
\address{Department of Mathematics, Tsuda University, 2-1-1 Tsuda-machi, Kodaira-shi, Tokyo 187-8577, Japan}
\email{ayminoue@tsuda.ac.jp}

\author{Naoki Kimura}
\address{Department of Mathematics, Graduate School of Fundamental Science and Engineering, Waseda University, 3-4-1 Okubo, Shinjuku-ku, Tokyo 169-8555, Japan}
\email{noverevitheuskyk@toki.waseda.jp}

\author{Ryo Nikkuni}
\address{Department of Mathematics, School of Arts and Sciences, Tokyo Woman's Christian University, 2-6-1 Zempukuji, Suginami-ku, Tokyo 167-8585, Japan}
\email{nick@lab.twcu.ac.jp}

\author{Kouki Taniyama}
\address{Department of Mathematics, School of Education, Waseda University, Nishi-Waseda 1-6-1, Shinjuku-ku, Tokyo, 169-8050, Japan}
\email{taniyama@waseda.jp}

\thanks{The first author was partially supported by Grant-in-Aid for Scientific Research (c) (No. 19K03476) , Japan Society for the Promotion of Science. The third author was partially supported by Grant-in-Aid for Scientific Research (c) (No. 19K03500) , Japan Society for the Promotion of Science. 
The fourth author was partially supported by Grant-in-Aid for Scientific Research (c) (No. 21K03260) and Grant-in-Aid for Scientific Research (A) (No. 21H04428), Japan Society for the Promotion of Science.}

\subjclass[2020]{Primary 05C10; Secondly 57K10.}

\date{}

\dedicatory{}

\keywords{crossing number, rotation number, plane immersed graph, Petersen graph, Heawood graph, Reduced Wu and generalized Simon invariant, Legendrian knot, Legendrian spatial graph, Thurston-Bennequin number}

\begin{abstract}

For any generic immersion of a Petersen graph into a plane, the number of crossing points between two edges of distance one is odd. The sum of the crossing numbers of all $5$-cycles is odd. The sum of the rotation numbers of all $5$-cycles is even. 
We show analogous results for $6$-cycles, $8$-cycles and $9$-cycles. 
For any Legendrian spatial embedding of a Petersen graph, there exists a $5$-cycle that is not an unknot with maximal Thurston-Bennequin number, and the sum of all Thurston-Bennequin numbers of the cycles is $7$ times the sum of all Thurston-Bennequin numbers of the $5$-cycles. We show analogous results for a Heawood graph. We also show some other results for some graphs. We characterize abstract graphs that has a generic immersion into a plane whose all cycles have rotation number $0$. 

\end{abstract}

\maketitle

\section{Introduction}\label{introduction} 

Let $G$ be a finite graph. We denote the set of all vertices of $G$ by $V(G)$ and the set of all edges by $E(G)$. 
We consider $G$ as a topological space in the usual way. Then a vertex of $G$ is a point of $G$ and an edge of $G$ is a subspace of $G$. 
A graph $G$ is said to be {\it simple} if it has no loops and no multiple edges. 
Suppose that $G$ has no multiple edges. Let $u$ and $v$ be mutually adjacent vertices of $G$. 
Then the edge of $G$ incident to both $u$ and $v$ is denoted by $uv$. 
Then $uv=vu$ as an unoriented edge. The orientation of $uv$ is given so that $u$ is the initial vertex and $v$ is the terminal vertex. 
Therefore $uv\neq vu$ as oriented edges. 
A {\it cycle} of $G$ is a subgraph of $G$ that is homeomorphic to a circle ${\mathbb S}^{1}$. 
A cycle with $k$ edges is said to be a {\it $k$-cycle}. 
The set of all cycles of $G$ is denoted by $\Gamma(G)$ and the set of all $k$-cycles of $G$ is denoted by $\Gamma_{k}(G)$. 

A {\it spatial embedding} of $G$ is an embedding of $G$ into the $3$-space ${\mathbb R}^{3}$. 
The image of a spatial embedding is said to be a {\it spatial graph}. 
The set of all spatial embeddings of $G$ is denoted by $SE(G)$. 
A {\it plane generic immersion} of $G$ is an immersion of $G$ into the plane ${\mathbb R}^{2}$ whose multiple points are only finitely many transversal double points between edges. 
Such a double point is said to be a {\it crossing point} or a {\it crossing}. 
The image of a plane generic immersion together with the distinction of the crossing points and the image of the degree $4$ vertices is said to be a {\it plane immersed graph}. 
Let $f:G\to{\mathbb R}^{2}$ be a plane generic immersion of $G$. Let $H$ be a subgraph of $G$. We denote the number of crossings of the restriction map $f|_{H}:H\to{\mathbb R}^{2}$ by $c(f(H))$. 
The set of all plane generic immersion of $G$ is denoted by $PGI(G)$. 

Let $K_{n}$ be a complete graph on $n$ vertices and $K_{m,n}$ a complete bipartite graph on $m+n$ vertices. 
It is shown in \cite{C-G} and \cite{Sachs} that for any spatial embedding $f:K_{6}\to {\mathbb R}^{3}$, the sum of all linking numbers of the links in $f(K_{6})$ is an odd number. 
Let $a_{2}(J)$ be the second coefficient of the Conway polynomial of a knot $J$. 
It is also shown in \cite{C-G} that for any spatial embedding $f:K_{7}\to {\mathbb R}^{3}$, the sum $\displaystyle{\sum_{\gamma\in\Gamma_{7}(K_{7})}a_{2}(f(\gamma))}$ is an odd number. 
See also \cite{Nikkuni} for refinements of these results, and \cite{S-S}\cite{S-S-S}\cite{Taniyama3} etc. for higher dimensional analogues. 
An analogous phenomenon appears in plane immersed graphs. 
A {\it self crossing} is a crossing of the same edge. An {\it adjacent crossing} is a crossing between two mutually adjacent edges. A {\it disjoint crossing} is a crossing between two mutually disjoint edges. 
It is known that for $G=K_{5}$ or $G=K_{3,3}$, the number of all disjoint crossings of a plane generic immersion of $G$ is always odd. See for example \cite[Proposition 2.1]{S-T} or \cite[Lemma 1.4.3]{Skopenkov}. Some theorems on plane immersed graphs are also stated in \cite{Skopenkov}. See also \cite{D-F-V} for related results. 
As analogous phenomenon we show the following results. 

Let $G$ be a finite graph with at least one cycle. The {\it girth} $g(G)$ of $G$ is the minimal lengths of the cycles of $G$. 
Namely every cycle of $G$ contains at least $g(G)$ edges and there is a $g(G)$-cycle of $G$. 
Let $G$ be a finite graph and $H,K$ connected subgraphs of $G$. The distance $d(H,K)$ of $H$ and $K$ in $G$ is defined to be the minimum number of edges of a path of $G$ joining $H$ and $K$. Then $d(H,K)=0$ if and only if $H\cup K$ is connected. Let $d$ and $e$ be mutually distinct edges of $G$. We note that $d(d,e)=0$ if and only if $d$ and $e$ are adjacent. Then $d(d,e)=1$ if and only if $d$ and $e$ are disjoint and there exists an edge $x$ of $G$ adjacent to both of them. If $g(G)\geq5$ then such $x$ is unique. 
Similarly $d(d,e)=2$ if and only if $d$ and $e$ are disjoint, no edge of $G$ is adjacent to both of them and there exist mutually adjacent edges $x$ and $y$ of $G$ such that $x$ is adjacent to $d$ and $y$ is adjacent to $e$. 
Let $k$ be a natural number. Let $D_{k}(G)$ be the set of all unordered pairs $(d,e)$ of edges of $G$ with $d(d,e)=k$.

\begin{Theorem}\label{theorem-K4}

Let $f:K_{4}\to {\mathbb R}^{2}$ be a plane generic immersion. 
Then 
\[
\sum_{\gamma\in\Gamma(K_{4})}c(f(\gamma))\equiv 0\pmod 2.
\]

\end{Theorem}

\begin{Theorem}\label{theorem-K33}

Let $f:K_{3,3}\to {\mathbb R}^{2}$ be a plane generic immersion. 
Then 
\[
\sum_{\gamma\in\Gamma_{4}(K_{3,3})}c(f(\gamma))\equiv\sum_{\gamma\in\Gamma_{6}(K_{3,3})}c(f(\gamma))\equiv 1\pmod 2.
\]

\end{Theorem}

\begin{figure}[htbp]
      \begin{center}
\scalebox{0.6}{\includegraphics*{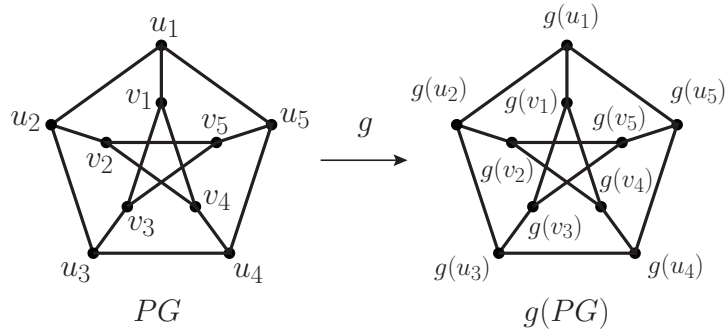}}
      \end{center}
   \caption{A plane generic immersion of $PG$}
  \label{PG}
\end{figure}

We denote a Petersen graph by $PG$. A Petersen graph $PG$ and a plane generic immersion $g:PG\to {\mathbb R}^{2}$ of $PG$ is illustrated in Figure \ref{PG}.

\begin{Theorem}\label{theorem-Petersen0}

Let $f:PG\to {\mathbb R}^{2}$ be a plane generic immersion. 
Then 
\[
\sum_{(d,e)\in D_{1}(PG)}|f(d)\cap f(e)|\equiv 1\pmod 2.
\]

\end{Theorem}

\begin{figure}[htbp]
      \begin{center}
\scalebox{0.6}{\includegraphics*{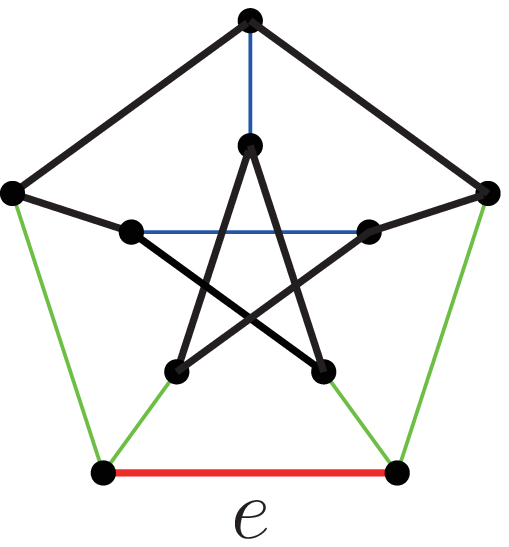}}
      \end{center}
   \caption{Distance $1$ edges form an $8$-cycle}
  \label{edge-distance-PG}
\end{figure}

Note that, for an edge $e$ of $PG$, the edges of $PG$ with distance $1$ with $e$ forms an $8$-cycle of $PG$, see Figure \ref{edge-distance-PG}.

\begin{figure}[htbp]
      \begin{center}
\scalebox{0.5}{\includegraphics*{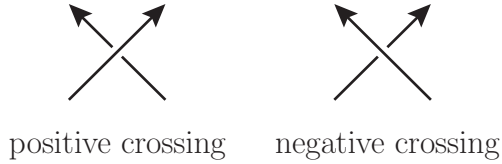}}
      \end{center}
   \caption{Positive crossing and negative crossing}
  \label{positive-negative-crossing}
\end{figure}

The modular equality in Theorem \ref{theorem-Petersen0} has an integral lift to an invariant of spatial embeddings of $PG$ as stated in Theorem \ref{theorem-Petersen1}. This is a kind of spatial graph invariants called {\it Reduced Wu and generalized Simon invariants}. See \cite{Taniyama1} \cite{Taniyama2} \cite{F-F-N}. 
We prepare some notions in order to state Theorem \ref{theorem-Petersen1}. 
Let $G$ be a finite graph and $f:G\to {\mathbb R}^{2}$ a plane generic immersion of $G$. Let $\pi:{\mathbb R}^{3}\to{\mathbb R}^{2}$ be a natural projection defined by $\pi(x,y,z)=(x,y)$. A spatial embedding $\varphi:G\to {\mathbb R}^{3}$ of $G$ is said to be a {\it lift} of $f$ if $f=\pi\circ\varphi$. The subset $f(G)$ of ${\mathbb R}^{2}$ together with the vertex information $f|_{V(G)}$ and over/under crossing information of $\varphi$ at each crossing of $f$ is said to be a {\it diagram} of $\varphi$ based on $f(G)$. 
Suppose $G$ is simple. Then a diagram of $\varphi$ restores $\varphi$ up to ambient isotopy of ${\mathbb R}^{3}$. 
Let $D$ be a diagram of $\varphi$. Suppose that each edge of $G$ is oriented. A crossing of $D$ is said to be a {\it positive crossing} or a {\it negative crossing} if it is as illustrated in Figure \ref{positive-negative-crossing}. 
Let $d$ and $e$ be mutually distinct edges of $G$. Let $c_{D}^{+}(d,e)$ be the number of positive crossings of $f(d)\cap f(e)$ and $c_{D}^{-}(d,e)$ the number of negative crossings of $f(d)\cap f(e)$. We set $\ell_{D}(d,e)=c_{D}^{+}(d,e)-c_{D}^{-}(d,e)$.

Let $u_{i}$ and $v_{i}$ be the vertices of a Petersen graph as illustrated in Figure \ref{PG} for $i=1,2,3,4,5$. We consider the suffixes modulo $5$. Namely $u_{-4}=u_{1}=u_{6}$, $u_{-3}=u_{2}=u_{7}$, $v_{-4}=v_{1}=v_{6}$ and so on. 
Then $E(PG)$ consists of $15$ edges $u_{i}u_{i+1}$, $u_{i}v_{i}$ and $v_{i}v_{i+2}$ for $i=1,2,3,4,5$. 
We consider that these edges are oriented. 
We define a map $\varepsilon:D_{1}(PG)\to{\mathbb Z}$ by 
$\varepsilon(u_{i}u_{i+1},u_{i+2}u_{i+3})=1$, 
$\varepsilon(u_{i}u_{i+1},u_{i-1}v_{i-1})=1$, 
$\varepsilon(u_{i}u_{i+1},u_{i+2}v_{i+2})=-1$, 
$\varepsilon(u_{i}u_{i+1},v_{j}v_{j+2})=1$,\\ 
$\varepsilon(u_{i}v_{i},u_{i\pm1}v_{i\pm1})=-1$, 
$\varepsilon(u_{i}v_{i},u_{i\pm2}v_{i\pm2})=1$, 
$\varepsilon(u_{i}v_{i},v_{i+1}v_{i+3})=-1$,\\ 
$\varepsilon(u_{i}v_{i},v_{i+2}v_{i+4})=1$ 
and 
$\varepsilon(v_{i}v_{i+2},v_{i+1}v_{i+3})=-1$ for $i,j\in\{1,2,3,4,5\}$. 
Let $\varphi:PG\to{\mathbb R}^{3}$ be a spatial embedding that is a lift of a plane generic immersion $f:PG\to{\mathbb R}^{2}$. 
Let $D$ be a diagram of $\varphi$ based on $f(PG)$. We set 
\[
{\mathcal L}(\varphi)=\sum_{(d,e)\in D_{1}(PG)}\varepsilon(d,e)\ell_{D}(d,e).
\]

\begin{Theorem}\label{theorem-Petersen1}

Let $\varphi:PG\to{\mathbb R}^{3}$ be a spatial embedding of a Petersen graph $PG$ that is a lift of a plane generic immersion $f:PG\to{\mathbb R}^{2}$ of $PG$. Let $D$ be a diagram of $\varphi$ based on $f(PG)$. Then ${\mathcal L}(\varphi)$ is a well-defined ambient isotopy invariant of $\varphi$ and we have 
\[
{\mathcal L}(\varphi)\equiv\sum_{(d,e)\in D_{1}(PG)}|f(d)\cap f(e)|\equiv 1\pmod 2.
\]

\end{Theorem}

\begin{Remark}\label{remark-Petersen1}
{\rm
Let $G$ be a finite graph and $k$ a natural number. 
Let $S_{k}(G)$ be a subcomplex of a $2$-dimensional complex $G\times G$ defined by 
\[
S_{k}(G)=\bigcup_{(d,e)\in D_{k}(G)}d\times e\cup e\times d.
\]
It is known that $S_{1}(K_{5})$ is homeomorphic to a closed orientable surface of genus $6$ and $S_{1}(K_{3,3})$ is homeomorphic to a closed orientable surface of genus $4$ \cite{Sarkaria}. 
Let $\varphi:G\to{\mathbb R}^{3}$ be a spatial embedding of $G$. 
Let $\tau_{\varphi}:S_{k}(G)\to{\mathbb S}^{2}$ be a Gauss map defined by
\[
\tau_{\varphi}(x,y)=\frac{\varphi(x)-\varphi(y)}{\|\varphi(x)-\varphi(y)\|}.
\]
It is known that the mapping degree of $\tau_{\varphi}:S_{1}(G)\to{\mathbb S}^{2}$ for $G=K_{5}$ or $K_{3,3}$ is equal to the {\it Simon invariant} of $\varphi$ up to sign \cite{Taniyama2}. 
By a straightforward consideration we see that $S_{1}(PG)$ is homeomorphic to a closed orientable surface of genus $16$, and for a spatial embedding $\varphi:PG\to{\mathbb R}^{3}$ we see that ${\mathcal L}(\varphi)$ is equal to the mapping degree of $\tau_{\varphi}:S_{1}(PG)\to{\mathbb S}^{2}$ up to sign. 
}
\end{Remark}

\begin{Theorem}\label{theorem-Petersen}

Let $f:PG\to {\mathbb R}^{2}$ be a plane generic immersion. 
Then 
\[
\sum_{\gamma\in\Gamma_{5}(PG)}c(f(\gamma))\equiv\sum_{\gamma\in\Gamma_{6}(PG)}c(f(\gamma))\equiv\sum_{\gamma\in\Gamma_{9}(PG)}c(f(\gamma))\equiv 1\pmod 2
\]
and 
\[
\sum_{\gamma\in\Gamma_{8}(PG)}c(f(\gamma))\equiv 0\pmod 4.
\]

\end{Theorem}

\begin{figure}[htbp]
      \begin{center}
\scalebox{0.6}{\includegraphics*{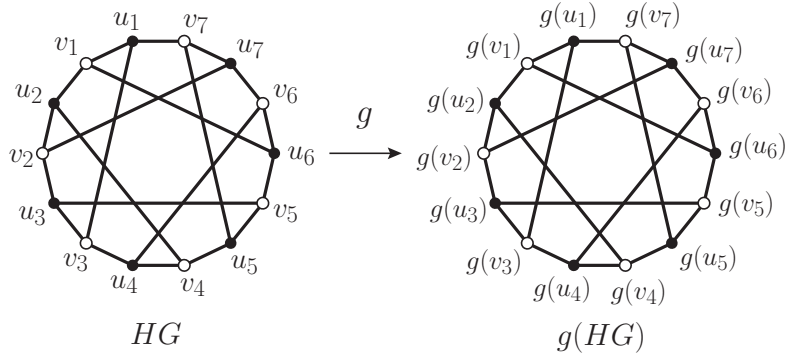}}
      \end{center}
   \caption{A plane generic immersion of $HG$}
  \label{HG}
\end{figure}

We denote a Heawood graph by $HG$. A Heawood graph $HG$ and a plane generic immersion $g:HG\to {\mathbb R}^{2}$ of $HG$ is illustrated in Figure \ref{HG}.

\begin{Theorem}\label{theorem-Heawood0}

Let $f:HG\to {\mathbb R}^{2}$ be a plane generic immersion. 
Then 
\[
\sum_{(d,e)\in D_{2}(HG)}|f(d)\cap f(e)|\equiv 1\pmod 2.
\]

\end{Theorem}

\begin{figure}[htbp]
      \begin{center}
\scalebox{0.6}{\includegraphics*{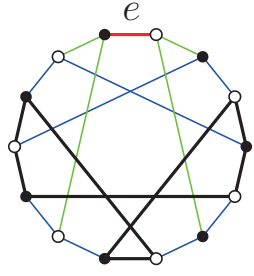}}
      \end{center}
   \caption{Distance $2$ edges form an $8$-cycle}
  \label{edge-distance-HG}
\end{figure}

Note that, for an edge $e$ of $HG$, the edges of $HG$ with distance $2$ with $e$ forms an $8$-cycle of $HG$, see Figure \ref{edge-distance-HG}.

We note that the modular equality in Theorem \ref{theorem-Heawood0} also has an integral lift to an invariant of spatial embeddings of $HG$. 
That is an invariant defined in \cite[Theorem 3.15]{F-F-N}. It is always an odd number \cite[Lemma 3.16]{F-F-N}.
The invariant is defined as follows. 
Let $u_{1},u_{2},u_{3},u_{4},u_{5},u_{6},u_{7}$ and $v_{1},v_{2},v_{3},v_{4},v_{5},v_{6},v_{7}$ be the vertices of a Heawood graph as illustrated in Figure \ref{HG}. We consider the suffixes modulo $7$. 
Then $E(HG)$ consists of $21$ edges $u_{i}v_{i}$, $u_{i}v_{i-1}$ and $v_{i}u_{i-2}$ for $i=1,2,3,4,5,6,7$. 
We consider that these edges are oriented. 
We define a map $\varepsilon:D_{1}(PG)\cup D_{2}(PG)\to{\mathbb Z}$ by 
$\varepsilon(u_{i}v_{i},u_{i\pm1}v_{i\pm1})=2$, 
$\varepsilon(u_{i}v_{i},u_{i\pm2}v_{i\pm2})=-2$, 
$\varepsilon(u_{i}v_{i},u_{i\pm3}v_{i\pm3})=-3$, 
$\varepsilon(u_{i}v_{i},u_{i+2}v_{i+1})=\varepsilon(u_{i}v_{i},u_{i-1}v_{i-2})=1$, \\
$\varepsilon(u_{i}v_{i},u_{i+3}v_{i+2})=\varepsilon(u_{i}v_{i},u_{i-2}v_{i-3})=2$, 
$\varepsilon(u_{i}v_{i},u_{i+4}v_{i+3})=5$, \\
$\varepsilon(u_{i}v_{i},v_{i+1}u_{i-1})=3$, 
$\varepsilon(u_{i}v_{i},v_{i+3}u_{i+1})=\varepsilon(u_{i}v_{i},v_{i-1}u_{i+4})=2$, \\
$\varepsilon(u_{i}v_{i},v_{i+4}u_{i+2})=\varepsilon(u_{i}v_{i},v_{i+5}u_{i+3})=-1$, \\
$\varepsilon(u_{i}v_{i-1},u_{i+1}v_{i})=\varepsilon(u_{i}v_{i-1},u_{i-1}v_{i-2})=2$, \\
$\varepsilon(u_{i}v_{i-1},u_{i+2}v_{i+1})=\varepsilon(u_{i}v_{i-1},u_{i-2}v_{i-3})=1$, \\
$\varepsilon(u_{i}v_{i-1},u_{i+3}v_{i+2})=\varepsilon(u_{i}v_{i-1},u_{i-3}v_{i-4})=-2$, \\
$\varepsilon(u_{i}v_{i-1},v_{i}u_{i-2})=\varepsilon(u_{i}v_{i-1},v_{i+1}u_{i-1})=2$, \\
$\varepsilon(u_{i}v_{i-1},v_{i+3}u_{i+1})=\varepsilon(u_{i}v_{i-1},v_{i+5}u_{i+3})=3$, \\
$\varepsilon(u_{i}v_{i-1},v_{i+4}u_{i+2})=3$, 
$\varepsilon(v_{i}u_{i-2},v_{i+1}u_{i-1})=\varepsilon(v_{i}u_{i-2},v_{i-1}u_{i-3})=5$, \\
$\varepsilon(v_{i}u_{i-2},v_{i+2}u_{i})=\varepsilon(v_{i}u_{i-2},v_{i-2}u_{i-4})=2$, 
and \\
$\varepsilon(v_{i}u_{i-2},v_{i+3}u_{i+1})=\varepsilon(v_{i}u_{i-2},v_{i-3}u_{i-5})=2$ 
for $i=1,2,3,4,5,6,7$. 
Let $\varphi:HG\to{\mathbb R}^{3}$ be a spatial embedding that is a lift of a plane generic immersion $f:HG\to{\mathbb R}^{2}$. 
Let $D$ be a diagram of $\varphi$ based on $f(HG)$. We set 
\[
{\mathcal L}(\varphi)=\sum_{(d,e)\in D_{1}(HG)\cup D_{2}(HG)}\varepsilon(d,e)\ell_{D}(d,e).
\]

\begin{Theorem}\label{theorem-Heawood1}

Let $\varphi:HG\to{\mathbb R}^{3}$ be a spatial embedding of a Heawood graph $HG$ that is a lift of a plane generic immersion $f:HG\to{\mathbb R}^{2}$ of $HG$. Let $D$ be a diagram of $\varphi$ based on $f(HG)$. Then we have 
\[
{\mathcal L}(\varphi)\equiv\sum_{(d,e)\in D_{2}(HG)}|f(d)\cap f(e)|\equiv 1\pmod 2.
\]

\end{Theorem}

\begin{Remark}\label{remark-Heawood1}
{\rm
The integral lift above involves disjoint edges of $HG$ with distance $1$. It is straightforward to check that there exists no integral lift involving only disjoint edges of $HG$ with distance $2$. As a related fact we see by a straightforward consideration that $S_{2}(HG)$ is homeomorphic to a closed non-orientable surface of non-orientable genus $30$, and for any spatial embedding $\varphi:HG\to{\mathbb R}^{3}$, the ${\mathbb Z}/2{\mathbb Z}$-mapping degree of $\tau_{\varphi}:S_{2}(HG)\to{\mathbb S}^{2}$ is equal to $1$. 
}
\end{Remark}

\begin{Theorem}\label{theorem-Heawood}

Let $f:HG\to {\mathbb R}^{2}$ be a plane generic immersion. 
Then 
\[
\sum_{\gamma\in\Gamma_{6}(HG)}c(f(\gamma))\equiv\sum_{\gamma\in\Gamma_{8}(HG)}c(f(\gamma))\equiv\sum_{\gamma\in\Gamma_{10}(HG)}c(f(\gamma))\equiv 1\pmod 2,
\]
\[
\sum_{\gamma\in\Gamma_{12}(HG)}c(f(\gamma))\equiv 0\pmod 4,
\]
and
\[
\sum_{\gamma\in\Gamma_{14}(HG)}c(f(\gamma))\equiv 0\pmod 2.
\]

\end{Theorem}

Let $f:G\to {\mathbb R}^{2}$ be a plane generic immersion of a finite graph $G$. Let $\gamma$ be a cycle of $G$. Suppose that $\gamma$ is given an orientation. Then $f(\gamma)$ is an oriented plane closed curve. We denote the rotation number of $f(\gamma)$ by ${\rm rot}(f(\gamma))$. 
We note that the parity of ${\rm rot}(f(\gamma))$ is independent of the choice of orientation of $\gamma$. 
It is easy to see that 
\[
{\rm rot}(f(\gamma))-c(f(\gamma))\equiv 1\pmod 2.
\]
We denote the number of elements of a finite set $X$ by $|X|$. 
The number of $k$-cycles of $G$ for $G=K_{4},K_{3,3}$ and $PG$ listed below is known and also easy to enumerate. 
The number of $k$-cycles of $HG$ listed below is also known. See for example \cite[4.2]{Sivaraman}. 
Since $|\Gamma(K_{4})|=7\equiv 1\pmod 2$, $|\Gamma_{4}(K_{3,3})|=9\equiv 1\pmod 2$, $|\Gamma_{6}(K_{3,3})|=6\equiv 0\pmod 2$, $|\Gamma_{5}(PG)|=12\equiv 0\pmod 2$, $|\Gamma_{6}(PG)|=10\equiv 0\pmod 2$, $|\Gamma_{8}(PG)|=15\equiv 1\pmod 2$, $|\Gamma_{9}(PG)|=20\equiv 0\pmod 2$, 
$|\Gamma_{6}(HG)|=28\equiv 0\pmod 2$, $|\Gamma_{8}(HG)|=21\equiv 1\pmod 2$, $|\Gamma_{10}(HG)|=84\equiv 0\pmod 2$, $|\Gamma_{12}(HG)|=56\equiv 0\pmod 2$ and $|\Gamma_{14}(HG)|=24\equiv 0\pmod 2$, we have the following immediate corollaries. 

\begin{Corollary}\label{corollary-K4}

Let $f:K_{4}\to {\mathbb R}^{2}$ be a plane generic immersion. 
Then 
\[
\sum_{\gamma\in\Gamma(K_{4})}{\rm rot}(f(\gamma))\equiv 1\pmod 2.
\]

\end{Corollary}

\begin{Corollary}\label{corollary-K33}

Let $f:K_{3,3}\to {\mathbb R}^{2}$ be a plane generic immersion. 
Then 
\[
\sum_{\gamma\in\Gamma_{4}(K_{3,3})}{\rm rot}(f(\gamma))\equiv 0\pmod 2
\]
and
\[
\sum_{\gamma\in\Gamma_{6}(K_{3,3})}{\rm rot}(f(\gamma))\equiv 1\pmod 2.
\]

\end{Corollary}

\begin{Corollary}\label{corollary-Petersen}

Let $f:PG\to {\mathbb R}^{2}$ be a plane generic immersion. 
Then 
\begin{align*}
&\sum_{\gamma\in\Gamma_{5}(PG)}{\rm rot}(f(\gamma))\equiv\sum_{\gamma\in\Gamma_{6}(PG)}{\rm rot}(f(\gamma))\equiv\sum_{\gamma\in\Gamma_{8}(PG)}{\rm rot}(f(\gamma))\\
\equiv
&\sum_{\gamma\in\Gamma_{9}(PG)}{\rm rot}(f(\gamma))\equiv 1\pmod 2.
\end{align*}

\end{Corollary}

\begin{Corollary}\label{corollary-Heawood}

Let $f:HG\to {\mathbb R}^{2}$ be a plane generic immersion. 
Then 
\[
\sum_{\gamma\in\Gamma_{6}(HG)}{\rm rot}(f(\gamma))\equiv\sum_{\gamma\in\Gamma_{10}(HG)}{\rm rot}(f(\gamma))\equiv 1\pmod 2
\]
and
\[
\sum_{\gamma\in\Gamma_{8}(HG)}{\rm rot}(f(\gamma))\equiv\sum_{\gamma\in\Gamma_{12}(HG)}{\rm rot}(f(\gamma))\equiv\sum_{\gamma\in\Gamma_{14}(HG)}{\rm rot}(f(\gamma))\equiv 0\pmod 2.
\]

\end{Corollary}

Furthermore we have the following theorem. 

\begin{Theorem}\label{theorem-zero-rotation}

Let $G$ be a finite graph. Then the following conditions are equivalent. 

\begin{enumerate}

\item[{\rm (1)}]

There is a plane generic immersion $f:G\to{\mathbb R}^{2}$ such that ${\rm rot}(f(\gamma))=0$ for every cycle $\gamma$ of $G$. 

\item[{\rm (2)}]

The graph $G$ does not have $K_{4}$ as a minor. 

\end{enumerate}

\end{Theorem}

Let $({\mathbb R}^{3},\xi_{\rm std})$ be the $3$-space with the standard contact structure. 
A {\it Legendrian knot} is a smooth knot in ${\mathbb R}^{3}$ that is tangent to the contact plane at each point. 
We consider Legendrian knots up to Legendrian isotopy. 
The Thurston-Bennequin number ${\rm tb}(K)$ and the rotation number ${\rm rot}(K)$ of a Legendrian knot $K$ are Legendrian isotopy invariants. 
A Legendrian knot $K$ is said to be a {\it trivial unknot} \cite{O-P} if it is a trivial knot as a classical knot and ${\rm tb}(K)=-1$. 
It is shown in \cite{O-P} that a finite graph $G$ has a Legendrian embedding with all cycles trivial unknots if and only if $G$ does not have $K_{4}$ as a minor. See also \cite{Tanaka} for related results. 
It is known that a trivial unknot has rotation number $0$. 
The rotation number of a Legendrian knot $K$ coincides with the rotation number of the plane immersed circle that is the image of $K$ under Lagrangian projection. 
Therefore the only if part is an immediate consequence of Corollary \ref{corollary-K4}. Namely we have the following corollary. 

\begin{Corollary}[\cite{O-P}]\label{corollary-K4-2}

Let $f:K_{4}\to {\mathbb R}^{3}$ be a Legendrian embedding. 
Then there is a cycle $\gamma$ of $K_{4}$ such that $f(\gamma)$ is not a trivial unknot. 

\end{Corollary}

Let $G$ be a finite graph and $f:G\to{\mathbb R}^{3}$ a Legendrian embedding. Then $f$ is said to be a {\it minimal embedding} \cite{O-P2} if $f(\gamma)$ is a trivial unknot for every $g(G)$-cycle $\gamma$ of $G$. Then we immediately have the following corollaries.

\begin{Corollary}\label{corollary-Petersen-Legendrian}

The Petersen graph $PG$ has no minimal Legendrian embedding. 

\end{Corollary}

\begin{Corollary}\label{corollary-Heawood-Legendrian}

The Heawood graph $HG$ has no minimal Legendrian embedding. 

\end{Corollary}

Let $G$ be a finite graph and $f:G\to{\mathbb R}^{3}$ a Legendrian embedding. The {\it total Thurston-Bennequin number} of $f$ is defined in \cite{O-P2} to be 
\[
TB(f)=\sum_{\gamma\in\Gamma(G)}tb(f(\gamma)).
\]
The following is also defined for a natural number $k$ in \cite{O-P2}. 
\[
TB_{k}(f)=\sum_{\gamma\in\Gamma_{k}(G)}tb(f(\gamma)).
\]
It is shown in \cite{O-P2} that $TB(f)$ is determined by $TB_{3}(f)$ when $G$ is a complete graph and by $TB_{4}(f)$ when $G$ is a complete bipartite graph. In this paper we extend these results to a Petersen graph and a Heawood graph.

\begin{Theorem}\label{theorem-Petersen-tb}

Let $f:PG\to{\mathbb R}^{3}$ be a Legendrian embedding. 
Then 
\[
TB_{6}(f)=TB_{5}(f),
\]
\[
TB_{8}(f)=2TB_{5}(f),
\]
and
\[
TB_{9}(f)=3TB_{5}(f).
\]
Therefore
\[
TB(f)=7TB_{5}(f).
\]

\end{Theorem}

\begin{Theorem}\label{theorem-Heawood-tb}

Let $f:HG\to{\mathbb R}^{3}$ be a Legendrian embedding. 
Then 
\[
TB_{8}(f)=TB_{6}(f),
\]
\[
TB_{10}(f)=5TB_{6}(f),
\]
\[
TB_{12}(f)=4TB_{6}(f),
\]
and
\[
TB_{14}(f)=2TB_{6}(f).
\]
Therefore
\[
TB(f)=13TB_{6}(f).
\]

\end{Theorem}

\section{Crossing numbers of cycles in a plane immersed graph}\label{crossing-numbers} 

\begin{Proposition}\label{proposition-crossing1}

Let $G$ be a finite graph and $\Lambda$ a set of subgraphs of $G$. Let $m$ be a positive integer. 
Suppose that the following {\rm (1)} and {\rm (2)} hold. 

\begin{enumerate}

\item[{\rm (1)}]

For every edge $e$ of $G$, the number of elements of $\Lambda$ containing $e$ is a multiple of $m$. 

\item[{\rm (2)}]

For every pair of edges $d$ and $e$ of $G$, the number of elements of $\Lambda$ containing both $d$ and $e$ is a multiple of $m$. 

\end{enumerate}

\noindent
Then for any plane generic immersion $f$ of $G$ 
\[
\sum_{\lambda\in\Lambda}c(f(\lambda))\equiv 0\pmod m.
\]

\end{Proposition}

\vskip 5mm

\noindent{\bf Proof.} Let $x$ be a self crossing of $f(G)$. Then by the condition (1) $x$ is counted a multiple of $m$ times in the sum. 
Let $x$ be an adjacent crossing or a disjoint crossing of $f(G)$. Then then by the condition (2) $x$ is counted a multiple of $m$ times in the sum. Therefore the sum is a multiple of $m$. 
$\Box$

\vskip 5mm

\noindent{\bf Proof of Theorem \ref{theorem-K4}.} We note that $\Gamma(K_{4})=\Gamma_{3}(K_{4})\cup\Gamma_{4}(K_{4})$, $|\Gamma_{3}(K_{4})|=4$ and $|\Gamma_{4}(K_{4})|=3$. Every edge of $K_{4}$ is contained in $2$ $3$-cycles and $2$ $4$-cycles. The total $4$ is a multiple of $2$. Every pair of mutually adjacent edges of $K_{4}$ is contained in a $3$-cycle and a $4$-cycle. The total $2$ is a multiple of $2$. Every pair of disjoint edges of $K_{4}$ is contained in no $3$-cycles and $2$ $4$-cycles. The total $2$ is a multiple of $2$. Then by Proposition \ref{proposition-crossing1} we have the result. 
$\Box$

\vskip 5mm

We note that the phenomenon described in Theorem \ref{theorem-K4} widely appears on graphs with certain symmetries. 
We show two of them below. The proofs are entirely analogous and we omit them.

\begin{Theorem}\label{theorem-K5}

Let $f:K_{5}\to {\mathbb R}^{2}$ be a plane generic immersion. 
Then 
\[
\sum_{\gamma\in\Gamma_{4}(K_{5})}c(f(\gamma))\equiv\sum_{\gamma\in\Gamma_{5}(K_{5})}c(f(\gamma))\equiv 0\pmod 2.
\]

\end{Theorem}

Let $m$ be a natural number. Let $T(m)$ be a graph of $3$ vertices and $3m$ edges such that each pair of vertices is joined by exactly $m$ multiple edges. 

\begin{Theorem}\label{theorem-multiple-triangle}
Let $m$ be a natural number. Let $f:T(m)\to{\mathbb R}^{2}$ be a plane generic immersion. Then 
\[
\sum_{\gamma\in\Gamma_{3}(T(m))}c(f(\gamma))\equiv 0\pmod m.
\]
\end{Theorem}

\vskip 5mm

\begin{figure}[htbp]
      \begin{center}
\scalebox{0.7}{\includegraphics*{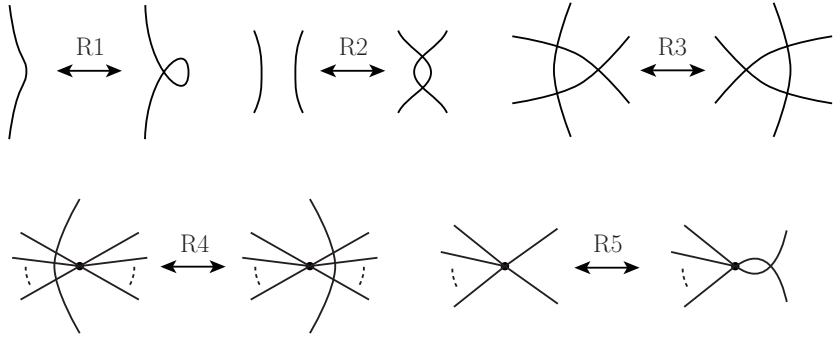}}
      \end{center}
   \caption{Reidemeister moves}
  \label{Reidemeister}
\end{figure}

It is known that any two plane generic immersions are transformed into each other up to self-homeomorphisms of $G$ and ${\mathbb R}^2$ by a finite sequence of local moves illustrated in Figure \ref{Reidemeister}. 
These moves are called {\it Reidemeister moves}. 

\begin{Proposition}\label{proposition-crossing2}

Let $G$ be a finite graph and $\Lambda$ a set of subgraphs of $G$. Let $m$ be a positive integer. 
The following conditions {\rm (A)} and {\rm (B)} are mutually equivalent. 

\begin{enumerate}

\item[{\rm (A)}]

For any plane generic immersions $f$ and $g$ of $G$, 
\[
\sum_{\lambda\in\Lambda}c(f(\lambda))\equiv \sum_{\lambda\in\Lambda}c(g(\lambda))\pmod m.
\]

\item[{\rm (B)}]

All of the following conditions {\rm (1)}, {\rm (2)}, {\rm (3)} and {\rm (4)} hold. 

\begin{enumerate}

\item[{\rm (1)}]

For every edge $e$ of $G$, the number of elements of $\Lambda$ containing $e$ is a multiple of $m$. 

\item[{\rm (2)}]

For every pair of edges $d$ and $e$ of $G$, twice the number of elements of $\Lambda$ containing both $d$ and $e$ is a multiple of $m$. 

\item[{\rm (3)}]

For every pair of a vertex $v$ and an edge $e$ of $G$, 
\[
\sum_{i=1}^{k}|\{\lambda\in\Lambda\mid\lambda\supset e\cup e_{i}\}|\equiv 0\pmod m,
\]
where $e_{1},\cdots,e_{k}$ are the edges of $G$ incident to $v$. 

\item[{\rm (4)}]

For every pair of mutually adjacent edges $d$ and $e$ of $G$, the number of elements of $\Lambda$ containing both $d$ and $e$ is a multiple of $m$. 

\end{enumerate}

\end{enumerate}

\end{Proposition}

\begin{Remark}\label{remark-crossing}
{\rm
If $\Lambda\subset\Gamma(G)$, then 
\[
\sum_{i=1}^{k}|\{\lambda\in\Lambda\mid\lambda\supset e\cup e_{i}\}|
=2\sum_{1\leq i<j\leq k}|\{\lambda\in\Lambda\mid\lambda\supset e\cup e_{i}\cup e_{j}\}|
\equiv 0\pmod 2.
\]
Therefore the condition (3) of (B) for $m=2$ automatically holds. 
}
\end{Remark}

\vskip 5mm

\begin{figure}[htbp]
      \begin{center}
\scalebox{0.4}{\includegraphics*{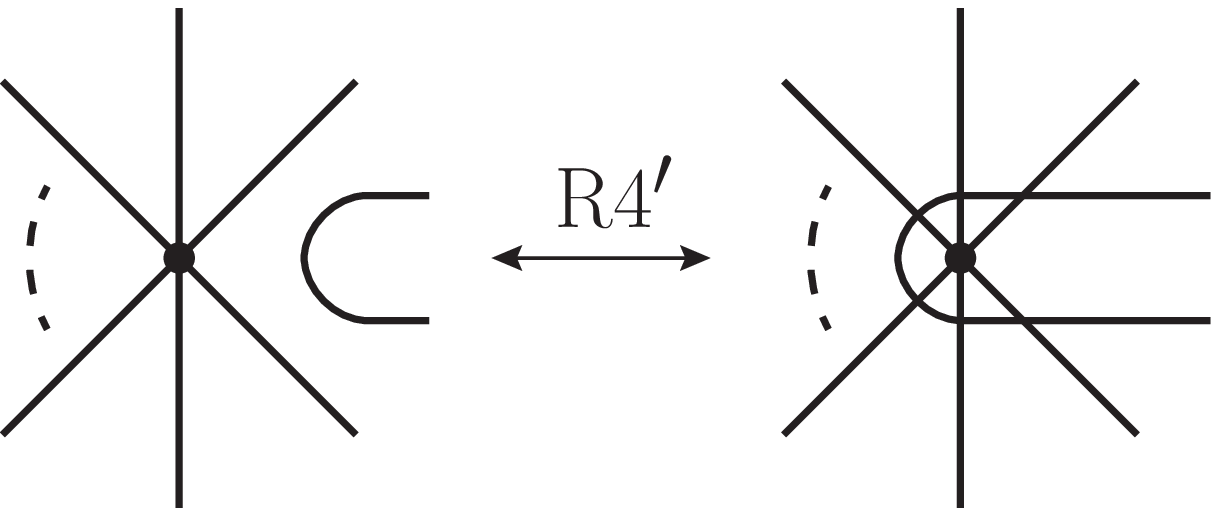}}
      \end{center}
   \caption{${\rm R4'}$}
  \label{Reidemeister2}
\end{figure}

\noindent{\bf Proof of Proposition \ref{proposition-crossing2}.} We set 
$\displaystyle{\tau_{\Lambda}(h)=\sum_{\lambda\in\Lambda}c(h(\lambda))}$ 
for a plane generic immersion $h$ of $G$. 
We note that a Reidemeister move ${\rm R4}$ in Figure \ref{Reidemeister} is realized by a Reidemeister move ${\rm R4'}$ in Figure \ref{Reidemeister2} and Reidemeister moves ${\rm R2}$ in Figure \ref{Reidemeister}. 
Therefore $f$ and $g$ are transformed into each other by ${\rm R1}$, ${\rm R2}$, ${\rm R3}$, ${\rm R4'}$ and ${\rm R5}$. 
A Reidemeister move ${\rm R1}$ create or annihilate a self crossing. Then we see that $\tau_{\Lambda}\pmod m$ is invariant under ${\rm R1}$ if and only if the condition (1) holds. 
A Reidemeister move ${\rm R2}$ create or annihilate two crossings. Suppose that they are self crossings, then $\tau_{\Lambda}\pmod m$ is invariant if the condition (1) holds. Suppose that they are both adjacent crossings or both disjoint crossings, then $\tau_{\Lambda}\pmod m$ is invariant if the condition (2) holds. 
If the condition (2) does not hold, then we can find $f$ and $g$ that differs by a Reidemeister move ${\rm R2}$ such that 
$\displaystyle{
\sum_{\lambda\in\Lambda}c(f(\lambda))-\sum_{\lambda\in\Lambda}c(g(\lambda))
}$
is not a multiple of $m$. 
A Reidemeister move ${\rm R3}$ does not change the number of self crossings, adjacent crossings and disjoint crossings for every edge, pair of mutually adjacent edges and pair of disjoint edges respectively. 
Therefore $\tau_{\Lambda}$ is always invariant under ${\rm R3}$. 
A Reidemeister move ${\rm R4'}$ create or annihilate crossings between an edge $e$ of $G$ and the edges of $G$ incident to $v$. 
Then we see that $\tau_{\Lambda}\pmod m$ is invariant under ${\rm R4'}$ if and only if the condition (3) holds. 
A Reidemeister move ${\rm R5}$ create or annihilate an adjacent crossing. Then we see that $\tau_{\Lambda}\pmod m$ is invariant under ${\rm R5}$ if and only if the condition (4) holds. This completes the proof. 
$\Box$

\vskip 5mm

\begin{figure}[htbp]
      \begin{center}
\scalebox{0.8}{\includegraphics*{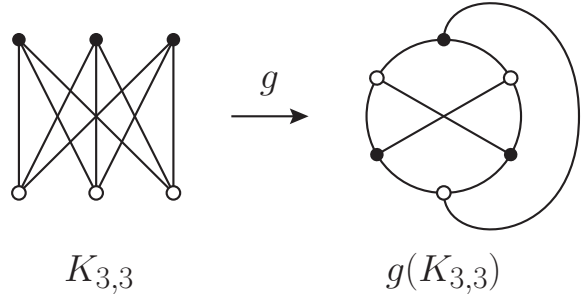}}
      \end{center}
   \caption{A plane generic immersion of $K_{3,3}$}
  \label{K33}
\end{figure}

\noindent{\bf Proof of Theorem \ref{theorem-K33}.} Let $g:K_{3,3}\to{\mathbb R}^{2}$ be a plane generic immersion illustrated in Figure \ref{K33}. We see that exactly one $4$-cycle of $K_{3,3}$ has a crossing under $g$ and exactly three $6$-cycles of $K_{3,3}$ has a crossing under $g$. Therefore 
\[
\sum_{\gamma\in\Gamma_{4}(K_{3,3})}c(g(\gamma))\equiv\sum_{\gamma\in\Gamma_{6}(K_{3,3})}c(g(\gamma))\equiv 1\pmod 2.
\]
Next we will check that the conditions (1), (2), (3) and (4) of Proposition \ref{proposition-crossing2} for $G=K_{3,3}$, $\Lambda=\Gamma_{4}(K_{3,3})$ or $\Lambda=\Gamma_{6}(K_{3,3})$, and $m=2$ hold. Each edge of $K_{3,3}$ is contained in exactly $4$ $4$-cycles and $4$ $6$-cycles of $K_{3,3}$. Therefore (1) holds. Since $m=2$, (2) automatically holds. 
Since $\Lambda\subset\Gamma(K_{3,3})$ and $m=2$, (3) holds by Remark \ref{remark-crossing}. 
Let $d$ and $e$ be mutually adjacent edges of $K_{3,3}$. Then we see that there exist exactly $2$ $4$-cycles and exactly $2$ $6$-cycles containing both of them. 
Therefore (4) holds. 
Then by Proposition \ref{proposition-crossing2} we have 
\[
\sum_{\gamma\in\Gamma_{4}(K_{3,3})}c(f(\gamma))\equiv \sum_{\gamma\in\Gamma_{4}(K_{3,3})}c(g(\gamma))\pmod 2
\]
and
\[
\sum_{\gamma\in\Gamma_{6}(K_{3,3})}c(f(\gamma))\equiv \sum_{\gamma\in\Gamma_{6}(K_{3,3})}c(g(\gamma))\pmod 2.
\]
Therefore we have 
\[
\sum_{\gamma\in\Gamma_{4}(K_{3,3})}c(f(\gamma))\equiv\sum_{\gamma\in\Gamma_{6}(K_{3,3})}c(f(\gamma))\equiv 1\pmod 2. 
\]
$\Box$

\vskip 5mm

\noindent{\bf Proof of Theorem \ref{theorem-Petersen0}.} We set 
\[
\kappa(f)=\sum_{(d,e)\in D_{1}(PG)}|f(d)\cap f(e)|
\]
for a plane generic immersion $f:PG\to{\mathbb R}^{2}$ of $PG$. 
Let $g:PG\to{\mathbb R}^{2}$ be a plane generic immersion illustrated in Figure \ref{PG}. We note that each crossing of $g(PG)$ is a disjoint crossing between two edges of $PG$ with distance $1$. Therefore we have 
\[
\kappa(g)=\sum_{(d,e)\in D_{1}(PG)}|g(d)\cap g(e)|=5\equiv 1\pmod 2.
\]
Next we will show that $\kappa(f) \pmod 2$ is invariant under Reidemeister moves. 
Since $\kappa$ do not count self crossings, it is invariant under ${\rm R1}$. The change of $\kappa$ under ${\rm R2}$ is $\pm2$ or $0$ and therefore $\kappa(f) \pmod 2$ is invariant under ${\rm R2}$. The third Reidemeister move ${\rm R3}$ do not change $\kappa$ itself. 
Let $e$ be an edge of $PG$ that is involved in a fourth Reidemeister move ${\rm R4}$. As we saw before, the edges of $PG$ with distance $1$ with $e$ forms an $8$-cycle of $PG$. Therefore we see that the change of $\kappa$ under ${\rm R4}$ is $\pm2$ or $0$ and therefore $\kappa(f) \pmod 2$ is invariant under ${\rm R4}$. Since $\kappa$ do not count adjacent crossings, it is invariant under ${\rm R5}$. Since $f$ and $g$ are transformed into each other by Reidemeister moves, we have 
\[
\kappa(f)\equiv\kappa(g)\equiv 1\pmod 2.
\]
$\Box$

\vskip 5mm

\noindent{\bf Proof of Theorem \ref{theorem-Petersen1}.} The proof that ${\mathcal L}(\varphi)$ is a well-defined ambient isotopy invariant taking its value in a odd number is entirely analogous to that of Simon invariants \cite{Taniyama1} and reduced Wu and generalized Simon invariants \cite{F-F-N}. The modular equality immediately follows from the definitions. 
$\Box$

\vskip 5mm

We prepare the following lemma for the proof of Theorem \ref{theorem-Petersen}.

\begin{Lemma}\label{lemma-PG}\

\begin{enumerate}

\item[{\rm (1)}]

For any vertices $x$ and $y$ of $PG$, there exists an isomorphism of $PG$ that maps $x$ to $y$. 

\item[{\rm (2)}]

For any edges $x$ and $y$ of $PG$, there exists an isomorphism of $PG$ that maps $x$ to $y$. 

\item[{\rm (3)}]

For any pairs of mutually adjacent edges $x,y$ and $z,w$ of $PG$, there exists an isomorphism of $PG$ that maps $x\cup y$ to $z\cup w$. 

\item[{\rm (4)}]

Let $x,y,z$ and $w$ be edges of $PG$ with $d(x,y)=d(z,w)=1$. 
Then there exists an isomorphism of $PG$ that maps $x\cup y$ to $z\cup w$. 

\item[{\rm (5)}]

Let $x,y,z$ and $w$ be edges of $PG$ with $d(x,y)=d(z,w)=2$. 
Then there exists an isomorphism of $PG$ that maps $x\cup y$ to $z\cup w$. 

\end{enumerate}

\end{Lemma}

\vskip 5mm

\noindent{\bf Proof.} Let $p:PG\to PG$ be an isomorphism defined by $p(u_{i})=u_{i+1}$ and $p(v_{i})=v_{i+1}$ for $i=1,2,3,4,5$. 
Let $q:PG\to PG$ be an isomorphism defined by $q(u_{i})=u_{5-i}$ and $q(v_{i})=v_{5-i}$ for $i=1,2,3,4,5$. 
Let $r:PG\to PG$ be an isomorphism defined by $r(u_{i})=v_{2i}$ and $r(v_{i})=u_{2i}$ for $i=1,2,3,4,5$. 
Let $s:PG\to PG$ be an isomorphism defined by $s(u_{1})=u_{1}$, $s(u_{2})=u_{2}$, $s(u_{5})=u_{5}$, $s(v_{1})=v_{1}$, $s(u_{3})=v_{2}$, $s(u_{4})=v_{5}$, $s(v_{2})=u_{3}$, $s(v_{5})=u_{4}$, $s(v_{3})=v_{4}$ and $s(v_{4})=v_{3}$. 
Then we see that all isomorphisms requested in (1) are generated by $p$ and $r$. 
We note that $q$ exchanges an edge $u_{5}u_{1}$ for an edge $u_{5}u_{4}$ and $s$ exchanges an edge $u_{2}v_{2}$ for an edge $u_{2}u_{3}$. 
Then by combining $p$ and $r$ we have all isomorphisms requested in (2) and (3). 
Suppose $d(x,y)=d(z,w)=1$. Then there is an edge $e$ of $PG$ adjacent to both $x$ and $y$. 
By (2) we map $e$ to an edge $u_{1}u_{2}$. We note that $q$ maps an edge $u_{2}u_{3}$ to $u_{3}u_{2}$. 
Then by combining other isomorphisms we have an isomorphism that maps $u_{1}u_{2}$ to $u_{2}u_{1}$. 
Then combining $s$ and these isomorphisms if necessary, $x\cup y$ is mapped to $u_{2}u_{3}\cup u_{5}u_{1}$. 
Similarly $z\cup w$ is mapped to $u_{2}u_{3}\cup u_{5}u_{1}$. This implies that $x\cup y$ is mapped to $z\cup w$ and (4) holds. 
Suppose $d(x,y)=d(z,w)=2$. Then there are mutually adjacent edges $d$ and $e$ of $PG$ such that $d$ is adjacent to $x$ and $e$ is adjacent to $y$. By (3) we map $d\cup e$ to $u_{1}u_{2}\cup u_{5}u_{1}$. Then we see that $x\cup y$ is mapped to $u_{2}v_{2}\cup u_{5}u_{4}$ or $u_{2}u_{3}\cup u_{5}v_{5}$. Similarly $z\cup w$ is mapped to $u_{2}v_{2}\cup u_{5}u_{4}$ or $u_{2}u_{3}\cup u_{5}v_{5}$. 
Since $s$ exchanges $u_{2}v_{2}\cup u_{5}u_{4}$ for $u_{2}u_{3}\cup u_{5}v_{5}$, we see that (5) holds. 
$\Box$

\vskip 5mm

Let $G$ be a finite graph and $k$ a natural number. Let $e$ be an edge of $G$. We set $\alpha_{k}(e,G)=|\{\gamma\in\Gamma_{k}(G)\mid\gamma\supset e\}|$. 
Let $d$ be another edge of $G$. 
We set $\alpha_{k}(d\cup e,G)=|\{\gamma\in\Gamma_{k}(G)\mid\gamma\supset d\cup e\}|$. 
Suppose that $d$ and $e$ are mutually disjoint and oriented. 
Let $\gamma$ be a cycle of $G$ containing $d\cup e$. 
The cycle $\gamma$ is said to be {\it coherent} with respect to $d$ and $e$ if the orientations of $d$ and $e$ are coherent in $\gamma$.  Otherwise $\gamma$ is said to be {\it incoherent} with respect to $d$ and $e$. 
Let $C_{k}(d\cup e,G)$ be the set of all cycles of $\Gamma_{k}(G)$ containing $d\cup e$ coherent with respect to $d$ and $e$. Let $I_{k}(d\cup e,G)$ be the set of all cycles of $\Gamma_{k}(G)$ containing $d\cup e$ incoherent with respect to $d$ and $e$. 
Then $\alpha_{k}(d\cup e,G)=|C_{k}(d\cup e,G)|+|I_{k}(d\cup e,G)|$. 
We set $\beta_{k}(d\cup e,G)=|C_{k}(d\cup e,G)|-|I_{k}(d\cup e,G)|$.

\begin{Lemma}\label{lemma-PG2}\

\begin{enumerate}

\item[{\rm (1)}]

Let $e$ be an edge of $PG$. Then $\alpha_{5}(e,PG)=4$, $\alpha_{6}(e,PG)=4$, $\alpha_{8}(e,PG)=8$ and $\alpha_{9}(e,PG)=12$. 

\item[{\rm (2)}]

Let $d$ and $e$ be mutually adjacent edges of $PG$. Then $\alpha_{5}(d\cup e,PG)=2$, $\alpha_{6}(d\cup e,PG)=2$, $\alpha_{8}(d\cup e,PG)=4$ and $\alpha_{9}(d\cup e,PG)=6$. 

\item[{\rm (3)}]

Let $d$ and $e$ be mutually disjoint oriented edges of $PG$ with $d(d,e)=1$. 
Then $\alpha_{5}(d\cup e,PG)=1$, $\alpha_{6}(d\cup e,PG)=1$, $\alpha_{8}(d\cup e,PG)=4$ and $\alpha_{9}(d\cup e,PG)=7$.
Suppose that the $5$-cycle of $PG$ containing $d\cup e$ is coherent with respect to $d$ and $e$. 
Then $\beta_{5}(d\cup e,PG)=1$, $\beta_{6}(d\cup e,PG)=1$, $\beta_{8}(d\cup e,PG)=2$ and $\beta_{9}(d\cup e,PG)=3$.

\item[{\rm (4)}]

Let $d$ and $e$ be mutually disjoint oriented edges of $PG$ with $d(d,e)=2$. 
Then $\alpha_{5}(d\cup e,PG)=0$, $\alpha_{6}(d\cup e,PG)=2$, $\alpha_{8}(d\cup e,PG)=4$, $\alpha_{9}(d\cup e,PG)=8$ and 
$\beta_{5}(d\cup e,PG)=\beta_{6}(d\cup e,PG)=\beta_{8}(d\cup e,PG)=\beta_{9}(d\cup e,PG)=0$.

\end{enumerate}

\end{Lemma}

\vskip 5mm

\noindent{\bf Proof.} There are $6$ pairs of mutually disjoint $5$-cycles of $PG$. 
They are \\
$(u_{1}u_{2}u_{3}u_{4}u_{5}u_{1}, v_{1}v_{3}v_{5}v_{2}v_{4}v_{1})$ and 
$(u_{i}u_{i+1}v_{i+1}v_{i+4}u_{i+4}u_{i}, v_{i}v_{i+2}u_{i+2}u_{i+3}v_{i+3}v_{i})$ \\
for $i=1,2,3,4,5$. Then $\Gamma_{5}(PG)$ consists of these $12$ $5$-cycles. 

The $10$ $6$-cycles $u_{i}u_{i+1}u_{i+2}u_{i+3}v_{i+3}v_{i}u_{i}$ and 
$u_{i}u_{i+1}v_{i+1}v_{i+4}v_{i+2}v_{i}u_{i}$\\
for $i=1,2,3,4,5$ are the elements of $\Gamma_{6}(PG)$. 

The $15$ $8$-cycles $u_{i}u_{i+1}u_{i+2}v_{i+2}v_{i}v_{i+3}u_{i+3}u_{i+4}u_{i}$, \\
$u_{i}u_{i+1}v_{i+1}v_{i+3}v_{i}v_{i+2}v_{i+4}u_{i+4}u_{i}$ and $u_{i}u_{i+1}v_{i+1}v_{i+3}u_{i+3}u_{i+2}v_{i+2}v_{i}u_{i}$\\ 
for $i=1,2,3,4,5$ are the elements of $\Gamma_{8}(PG)$. 

The $20$ $9$-cycles \\
$u_{i}u_{i+1}v_{i+1}v_{i+3}u_{i+3}u_{i+2}v_{i+2}v_{i+4}u_{i+4}u_{i}$, 
$u_{i}u_{i+1}v_{i+1}v_{i+4}v_{i+2}u_{i+2}u_{i+3}v_{i+3}v_{i}u_{i}$, \\
$u_{i}u_{i+1}u_{i+2}v_{i+2}v_{i+4}v_{i+1}v_{i+3}u_{i+3}u_{i+4}u_{i}$ and 
$u_{i}u_{i+1}u_{i+2}u_{i+3}v_{i+3}v_{i+1}v_{i+4}v_{i+2}v_{i}u_{i}$ \\
for $i=1,2,3,4,5$ are the elements of $\Gamma_{9}(PG)$. 

Then by Lemma \ref{lemma-PG} (2) we see that counting cycles for a particular edge $e$ of $PG$ will show (1). Similarly we have (2), (3) and (4) by Lemma \ref{lemma-PG} (3), (4) and (5) respectively. 
$\Box$

\vskip 5mm

Summarizing the statements in Lemma \ref{lemma-PG2} we have Table \ref{table-Petersen-cycles}.

\begin{table}[htb]
\centering
\caption{number of $k$-cycles of $PG$}
\scalebox{0.75}
{
\begin{tabular}{|c|c|c|c|c|c|c|c|c|}  \hline
\multirow{2}{*}{$k$} & 
\multirow{2}{*}{$|\Gamma_{k}(PG)|$} & 
\multirow{2}{*}{$|\Gamma_{k}(PG)|\cdot k$} & 
\multirow{2}{*}{$\alpha_{k}(e,PG)$} & 
\multicolumn{3}{c}{$\alpha_{k}(d\cup e,PG)$} & 
\multicolumn{2}{|c|}{$\beta_{k}(d\cup e,PG)$} 
\\ \cline{5-9}
& & & & $d(d,e)=0$ & $d(d,e)=1$ & $d(d,e)=2$ &
$d(d,e)=1$ & $d(d,e)=2$ 

\\ \hline
    $5$ & $12$ & $60$ & $4$ & $2$ & $1$ & $0$ & $1$ & $0$ \\ \hline
    $6$ & $10$ & $60$ & $4$ & $2$ & $1$ & $2$ & $1$ & $0$ \\ \hline
    $8$ & $15$ & $120$ & $8$ & $4$ & $4$ & $4$ & $2$ & $0$ \\ \hline
    $9$ & $20$ & $180$ & $12$ & $6$ & $7$ & $8$ & $3$ & $0$ \\ \hline
\end{tabular}
}
  \label{table-Petersen-cycles}
\end{table}

\vskip 5mm

\noindent{\bf Proof of Theorem \ref{theorem-Petersen}.} Let $g:PG\to{\mathbb R}^{2}$ be a plane generic immersion illustrated in Figure \ref{PG}. Since each crossing of $g(PG)$ is a crossing between two edges of $PG$ with distance $1$ and $\alpha_{5}(d\cup e,PG)=1$ for edges $d$ and $e$ of $PG$ with $d(d,e)=1$, we have 
\[
\sum_{\gamma\in\Gamma_{5}(PG)}c(g(\gamma))=5\equiv 1\pmod 2.
\]
Since $\alpha_{6}(d\cup e,PG)=1$ and $\alpha_{9}(d\cup e,PG)=7$ for edges $d$ and $e$ of $PG$ with $d(d,e)=1$, we have
\[
\sum_{\gamma\in\Gamma_{6}(PG)}c(g(\gamma))=5\equiv 1\pmod 2
\]
and
\[
\sum_{\gamma\in\Gamma_{9}(PG)}c(g(\gamma))=35\equiv 1\pmod 2.
\]
Next we will check that the conditions (1), (2), (3) and (4) of Proposition \ref{proposition-crossing2} for $G=PG$, $\Lambda=\Gamma_{5}(PG)$, $\Gamma_{6}(PG)$ or $\Gamma_{9}(PG)$ and $m=2$ hold. 
Since $\alpha_{5}(e,PG)=\alpha_{6}(e,PG)=4$ and $\alpha_{9}(e,PG)=12$ for every edge $e$ of $PG$, we see that (1) holds. 
For $m=2$, (2) automatically holds. 
Since $\Lambda\subset\Gamma(PG)$ (3) holds for $m=2$ by Remark \ref{remark-crossing}. 
Let $d$ and $e$ be mutually adjacent edges of $PG$. 
Then $\alpha_{5}(d\cup e,PG)=\alpha_{6}(d\cup e,PG)=2$ and $\alpha_{9}(d\cup e,PG)=6$. 
Therefore (4) holds. 
Then by Proposition \ref{proposition-crossing2} we have 
\[
\sum_{\gamma\in\Gamma_{5}(PG)}c(f(\gamma))\equiv \sum_{\gamma\in\Gamma_{5}(PG)}c(g(\gamma))\pmod 2,
\]
\[
\sum_{\gamma\in\Gamma_{6}(PG)}c(f(\gamma))\equiv \sum_{\gamma\in\Gamma_{6}(PG)}c(g(\gamma))\pmod 2
\]
and
\[
\sum_{\gamma\in\Gamma_{9}(PG)}c(f(\gamma))\equiv \sum_{\gamma\in\Gamma_{9}(PG)}c(g(\gamma))\pmod 2.
\]
Therefore we have 
\[
\sum_{\gamma\in\Gamma_{5}(PG)}c(f(\gamma))\equiv\sum_{\gamma\in\Gamma_{6}(PG)}c(f(\gamma))\equiv\sum_{\gamma\in\Gamma_{9}(PG)}c(f(\gamma))\equiv 1\pmod 2. 
\]
We note that $\alpha_{8}(e,PG)=8$ for every edge $e$ of $PG$ and $\alpha_{8}(d\cup e,PG)=4$ for any distinct edges $d$ and $e$. 
Then by Proposition \ref{proposition-crossing1} we have 
\[
\sum_{\gamma\in\Gamma_{8}(PG)}c(f(\gamma))\equiv 0\pmod 4.
\]
$\Box$

\vskip 5mm

\begin{figure}[htbp]
      \begin{center}
\scalebox{0.6}{\includegraphics*{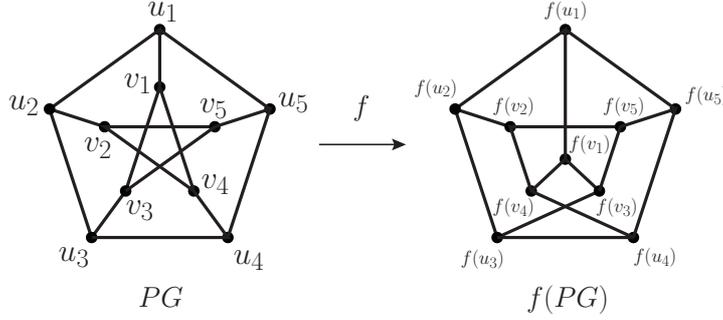}}
      \end{center}
   \caption{Another plane generic immersion of $PG$}
  \label{PG2}
\end{figure}

\begin{Example}\label{example-PG}
{\rm
Another plane generic immersion $f:PG\to {\mathbb R}^{2}$ of $PG$ is illustrated in Figure \ref{PG2}. The crossing number $c(f(PG))=2$ is known to be minimal among all plane generic immersions of $PG$. Note that the upper crossing of $f(PG)$ in Figure \ref{PG2} is a crossing between distance $2$ edges of $PG$ and the lower crossing is a crossing between distance $1$ edges of $PG$. 
Then we see that $f$ satisfies all modular equalities in Theorem \ref{theorem-Petersen0} and Theorem \ref{theorem-Petersen}. 
}
\end{Example}

\vskip 5mm

\noindent{\bf Proof of Theorem \ref{theorem-Heawood0}.} We set 
\[
\kappa(f)=\sum_{(d,e)\in D_{2}(HG)}|f(d)\cap f(e)|
\]
for a plane generic immersion $f:HG\to{\mathbb R}^{2}$ of $HG$. 
Let $g:HG\to{\mathbb R}^{2}$ be a plane generic immersion illustrated in Figure \ref{HG}. We note that $g(HG)$ has $7$ crossings between distance $1$ edges and $7$ crossings between distance $2$ edges. Therefore we have 
\[
\kappa(g)=\sum_{(d,e)\in D_{2}(HG)}|g(d)\cap g(e)|=7\equiv 1\pmod 2.
\]
Next we will show that $\kappa(f) \pmod 2$ is invariant under Reidemeister moves. 
Since $\kappa$ do not count self crossings, it is invariant under ${\rm R1}$. The change of $\kappa$ under ${\rm R2}$ is $\pm2$ or $0$ and therefore $\kappa(f) \pmod 2$ is invariant under ${\rm R2}$. The third Reidemeister move ${\rm R3}$ do not change $\kappa$ itself. 
Let $e$ be an edge of $HG$ that is involved in a fourth Reidemeister move ${\rm R4}$. As we saw before, the edges of $HG$ with distance $2$ with $e$ forms an $8$-cycle of $HG$. Therefore we see that the change of $\kappa$ under ${\rm R4}$ is $\pm2$ or $0$ and therefore $\kappa(f) \pmod 2$ is invariant under ${\rm R4}$. Since $\kappa$ do not count adjacent crossings, it is invariant under ${\rm R5}$. Since $f$ and $g$ are transformed into each other by Reidemeister moves, we have 
\[
\kappa(f)\equiv\kappa(g)\equiv 1\pmod 2.
\]
$\Box$

\vskip 5mm

\noindent{\bf Proof of Theorem \ref{theorem-Heawood1}.} We note that 
\[
\ell_{D}(d,e)=c_{D}^{+}(d,e)-c_{D}^{-}(d,e)\equiv c_{D}^{+}(d,e)+c_{D}^{-}(d,e)\equiv|d\cap e|\pmod 2.
\]
We also note that $\varepsilon(d,e)$ is an odd number if $d(d,e)=1$ and it is an even number if $d(d,e)=2$. 
Therefore we have the modular equality. 
$\Box$

\begin{Lemma}\label{lemma-HG}\

\begin{enumerate}

\item[{\rm (1)}]

For any vertices $x$ and $y$ of $HG$, there exists an isomorphism of $HG$ that maps $x$ to $y$. 

\item[{\rm (2)}]

For any edges $x$ and $y$ of $HG$, there exists an isomorphism of $HG$ that maps $x$ to $y$. 

\item[{\rm (3)}]

For any pairs of mutually adjacent edges $x,y$ and $z,w$ of $HG$, there exists an isomorphism of $HG$ that maps $x\cup y$ to $z\cup w$. 

\item[{\rm (4)}]

Let $x,y,z$ and $w$ be edges of $HG$ with $d(x,y)=d(z,w)=1$. 
Then there exists an isomorphism of $HG$ that maps $x\cup y$ to $z\cup w$. 

\item[{\rm (5)}]

Let $x,y,z$ and $w$ be edges of $HG$ with $d(x,y)=d(z,w)=2$. 
Then there exists an isomorphism of $HG$ that maps $x\cup y$ to $z\cup w$. 

\end{enumerate}

\end{Lemma}

\vskip 5mm

\noindent{\bf Proof.} A Heawood graph has the following symmetries. 
It is $7$-periodic with respect to the isomorphism $\alpha:HG\to HG$ defined by $\alpha(u_{i})=u_{i+1}$ and $\alpha(v_{i})=v_{i+1}$. 
It is $2$-periodic with respect to the isomorphism $\beta:HG\to HG$ defined by $\beta(u_{i})=v_{8-i}$ and $\beta(v_{i})=u_{8-i}$. 
It is $3$-periodic with $2$ fixed vertices with respect to the isomorphism $\gamma:HG\to HG$ defined by 
$\gamma(u_{1})=u_{1}$, $\gamma(v_{6})=v_{6}$, 
$\gamma(v_{1})=v_{7}$, $\gamma(v_{7})=v_{3}$, $\gamma(v_{3})=v_{1}$, 
$\gamma(u_{2})=u_{5}$, $\gamma(u_{5})=u_{3}$, $\gamma(u_{3})=u_{2}$, 
$\gamma(v_{2})=v_{4}$, $\gamma(v_{4})=v_{5}$, $\gamma(v_{5})=v_{2}$ 
and 
$\gamma(u_{4})=u_{6}$, $\gamma(u_{6})=u_{7}$, $\gamma(u_{7})=u_{4}$. 
Then $\alpha$ and $\beta$ generate all isomorphisms requested in (1). 
We note that $\gamma$ cyclically maps the edges incident to $u_{1}$. 
Then we see that these isomorphisms generate all isomorphisms requested in (2) and (3). 
Let $\varphi:HG\to HG$ be an isomorphism defined by 
$\varphi(u_{1})=u_{1}$, $\varphi(v_{7})=v_{7}$, 
$\varphi(u_{7})=u_{7}$, $\varphi(v_{4})=v_{4}$, 
$\varphi(u_{5})=u_{5}$, $\varphi(v_{5})=v_{5}$, 
$\varphi(v_{1})=v_{3}$, $\varphi(v_{3})=v_{1}$, 
$\varphi(u_{2})=u_{4}$, $\varphi(u_{4})=u_{2}$, 
$\varphi(v_{2})=v_{6}$, $\varphi(v_{6})=v_{2}$ 
and 
$\varphi(u_{3})=u_{6}$, $\varphi(u_{6})=u_{3}$. 
We note that $\varphi$ fixes the edges $u_{1}v_{7}$ and $u_{7}v_{7}$ and exchanges $u_{1}v_{1}$ for $u_{1}v_{3}$. 
Suppose $d(x,y)=d(z,w)=1$. Let $e$ be an edge of $HG$ adjacent to both of $x$ and $y$. By (2) we map $e$ to $u_{1}v_{7}$. Then by $\varphi$ and $\beta$ we map $x\cup y$ to $u_{1}v_{1}\cup u_{7}v_{7}$. Similarly $z\cup w$ is mapped to $u_{1}v_{1}\cup u_{7}v_{7}$. Thus we have an isomorphism requested in (4). 
Let $\psi:HG\to HG$ be an isomorphism defined by 
$\psi(v_{1})=v_{1}$, $\psi(u_{1})=u_{1}$, 
$\psi(v_{7})=v_{7}$, $\psi(u_{7})=u_{7}$, 
$\psi(v_{3})=v_{3}$, $\psi(u_{5})=u_{5}$, 
$\psi(u_{2})=u_{6}$, $\psi(u_{6})=u_{2}$, 
$\psi(v_{2})=v_{6}$, $\psi(v_{6})=v_{2}$, 
$\psi(u_{3})=u_{4}$, $\psi(u_{4})=u_{3}$, 
and 
$\psi(v_{4})=v_{5}$, $\psi(v_{5})=v_{4}$. 
We note that $\psi$ fixes the path $v_{1}u_{1}v_{7}u_{7}$ and exchanges $v_{1}u_{2}$ for $v_{1}u_{6}$. 
Suppose $d(x,y)=d(z,w)=2$. Let $d$ and $e$ be mutually adjacent edges of $HG$ such that $x\cup d\cup e\cup y$ is a path. 
We map $d\cup e$ to $u_{1}v_{7}\cup u_{7}v_{7}$ by (3). 
Then by $\varphi$ we map $x$ or $y$ to $u_{1}v_{1}$. Then by $\beta$ and $\psi$ the path $x\cup d\cup e\cup y$ is mapped to $u_{2}v_{1}u_{1}v_{7}u_{7}$. Thus $x\cup y$ is mapped to $u_{2}v_{1}\cup v_{7}u_{7}$. 
Similarly $z\cup w$ is mapped to $u_{2}v_{1}\cup v_{7}u_{7}$. Thus we have an isomorphism requested in (5). 
$\Box$

\begin{Lemma}\label{lemma-HG2}\

\begin{enumerate}

\item[{\rm (1)}]

Let $e$ be an edge of $HG$. Then $\alpha_{6}(e,HG)=8$, $\alpha_{8}(e,HG)=8$, $\alpha_{10}(e,HG)=40$, $\alpha_{12}(e,HG)=32$ and $\alpha_{14}(e,HG)=16$. 

\item[{\rm (2)}]

Let $d$ and $e$ be mutually adjacent edges of $HG$. Then $\alpha_{6}(d\cup e,HG)=4$, $\alpha_{8}(d\cup e,HG)=4$, $\alpha_{10}(d\cup e,HG)=20$, $\alpha_{12}(d\cup e,HG)=16$ and $\alpha_{14}(d\cup e,HG)=8$. 

\item[{\rm (3)}]

Let $d$ and $e$ be mutually disjoint oriented edges of $HG$ with $d(d,e)=1$. 
Then $\alpha_{6}(d\cup e,HG)=2$, $\alpha_{8}(d\cup e,HG)=2$, $\alpha_{10}(d\cup e,HG)=18$, $\alpha_{12}(d\cup e,HG)=16$ and $\alpha_{14}(d\cup e,HG)=12$.
Suppose that the $6$-cycle of $HG$ containing $d\cup e$ is coherent with respect to $d$ and $e$. 
Then $\beta_{6}(d\cup e,HG)=2$, $\beta_{8}(d\cup e,HG)=2$, $\beta_{10}(d\cup e,HG)=10$, $\beta_{12}(d\cup e,HG)=8$ and $\beta_{14}(d\cup e,HG)=4$.

\item[{\rm (4)}]

Let $d$ and $e$ be mutually disjoint edges of $HG$ with $d(d,e)=2$. 
Then $\alpha_{5}(d\cup e,HG)=0$, $\alpha_{6}(d\cup e,HG)=2$, $\alpha_{8}(d\cup e,HG)=4$, $\alpha_{8}(d\cup e,HG)=4$ and $\alpha_{9}(d\cup e,HG)=8$.
Suppose that the $6$-cycle of $HG$ containing $d\cup e$ is coherent with respect to $d$ and $e$. 
Then $\beta_{6}(d\cup e,HG)=1$, $\beta_{8}(d\cup e,HG)=1$, $\beta_{10}(d\cup e,HG)=5$, $\beta_{12}(d\cup e,HG)=4$ and $\beta_{14}(d\cup e,HG)=2$.

\end{enumerate}

\end{Lemma}

\vskip 5mm

\noindent{\bf Proof.} We note that $HG$ is a bipartite graph on $14$ vertices and $g(HG)=6$. Therefore 
\[
\Gamma(HG)=\bigcup_{k=3}^{7}\Gamma_{2k}(HG).
\]

There are $28$ $6$-cycles of $HG$. 
They are \\
$u_{i}v_{i}u_{i+1}v_{i+1}u_{i+2}v_{i+2}u_{i}$, $u_{i}v_{i}u_{i+1}v_{i+3}u_{i+3}v_{i+2}u_{i}$, $u_{i}v_{i+2}u_{i+3}v_{i+3}u_{i+4}v_{i+6}u_{i}$ \\
and $u_{i}v_{i+2}u_{i+2}v_{i+4}u_{i+4}v_{i+6}u_{i}$ 
for $i=1,2,3,4,5,6,7$. 

There are $21$ $8$-cycles of $HG$. 
They are \\
$u_{i+6}v_{i+6}u_{i}v_{i}u_{i+5}v_{i+4}u_{i+2}v_{i+1}u_{i+6}$, 
$v_{i+5}u_{i+6}v_{i+6}u_{i}v_{i}u_{i+1}v_{i+3}u_{i+3}v_{i+5}$ and\\ 
$u_{i+1}v_{i+1}u_{i+2}v_{i+4}u_{i+5}v_{i+5}u_{i+3}v_{i+3}u_{i+1}$ 
for $i=1,2,3,4,5,6,7$. 

There are $84$ $10$-cycles of $HG$. 
They are \\
$v_{i+4}u_{i+5}v_{i+5}u_{i+6}v_{i+6}u_{i}v_{i}u_{i+1}v_{i+1}u_{i+2}v_{i+4}$, \\
$u_{i+1}v_{i+1}u_{i+2}v_{i+2}u_{i+3}v_{i+5}u_{i+5}v_{i+4}u_{i+4}v_{i+3}u_{i+1}$, \\
$u_{i}v_{i+2}u_{i+2}v_{i+1}u_{i+6}v_{i+5}u_{i+5}v_{i+4}u_{i+4}v_{i+6}u_{i}$, \\
$u_{i}v_{i+2}u_{i+2}v_{i+1}u_{i+1}v_{i}u_{i+5}v_{i+4}u_{i+4}v_{i+6}u_{i}$, \\
$u_{i}v_{i}u_{i+1}v_{i+1}u_{i+6}v_{i+6}u_{i+4}v_{i+3}u_{i+3}v_{i+2}u_{i}$, \\
$u_{i}v_{i}u_{i+5}v_{i+5}u_{i+6}v_{i+6}u_{i+4}v_{i+3}u_{i+3}v_{i+2}u_{i}$, \\
$u_{i}v_{i}u_{i+1}v_{i+3}u_{i+4}v_{i+6}u_{i+6}v_{i+5}u_{i+3}v_{i+2}u_{i}$, \\
$u_{i+6}v_{i+6}u_{i}v_{i}u_{i+5}v_{i+5}u_{i+3}v_{i+3}u_{i+1}v_{i+1}u_{i+6}$, \\
$u_{i}v_{i+2}u_{i+3}v_{i+3}u_{i+1}v_{i+1}u_{i+2}v_{i+4}u_{i+4}v_{i+6}u_{i}$, \\
$u_{i}v_{i+2}u_{i+2}v_{i+4}u_{i+5}v_{i+5}u_{i+3}v_{i+3}u_{i+4}v_{i+6}u_{i}$, \\
$u_{i}v_{i}u_{i+5}v_{i+4}u_{i+4}v_{i+6}u_{i+6}v_{i+1}u_{i+2}v_{i+2}u_{i}$ and\\
$v_{i}u_{i+1}v_{i+3}u_{i+3}v_{i+5}u_{i+6}v_{i+1}u_{i+2}v_{i+4}u_{i+5}v_{i}$ 
for $i=1,2,3,4,5,6,7$. 

There are $56$ $12$-cycles of $HG$. 
They are \\
$u_{i+6}v_{i+6}u_{i}v_{i}u_{i+5}v_{i+4}u_{i+4}v_{i+3}u_{i+3}v_{i+2}u_{i+2}v_{i+1}u_{i+6}$, \\
$v_{i+5}u_{i+6}v_{i+6}u_{i}v_{i}u_{i+1}v_{i+3}u_{i+4}v_{i+4}u_{i+2}v_{i+2}u_{i+3}v_{i+5}$, \\
$u_{i}v_{i+2}u_{i+2}v_{i+1}u_{i+1}v_{i+3}u_{i+3}v_{i+5}u_{i+5}v_{i+4}u_{i+4}v_{i+6}u_{i}$, \\
$u_{i}v_{i}u_{i+1}v_{i+1}u_{i+2}v_{i+4}u_{i+4}v_{i+6}u_{i+6}v_{i+5}u_{i+3}v_{i+2}u_{i}$, \\
$u_{i}v_{i}u_{i+1}v_{i+3}u_{i+4}v_{i+6}u_{i+6}v_{i+5}u_{i+5}v_{i+4}u_{i+2}v_{i+2}u_{i}$, \\
$u_{i}v_{i}u_{i+1}v_{i+3}u_{i+3}v_{i+5}u_{i+6}v_{i+6}u_{i+4}v_{i+4}u_{i+2}v_{i+2}u_{i}$, \\
$u_{i}v_{i}u_{i+5}v_{i+4}u_{i+2}v_{i+1}u_{i+6}v_{i+6}u_{i+4}v_{i+3}u_{i+3}v_{i+2}u_{i}$ and\\
$u_{i}v_{i+2}u_{i+3}v_{i+5}u_{i+5}v_{i+4}u_{i+2}v_{i+1}u_{i+1}v_{i+3}u_{i+4}v_{i+6}u_{i}$ \\
for $i=1,2,3,4,5,6,7$. 

There are $24$ $14$-cycles of $HG$. \\
They are $u_{1}v_{1}u_{2}v_{2}u_{3}v_{3}u_{4}v_{4}u_{5}v_{5}u_{6}v_{6}u_{7}v_{7}u_{1}$, 
$u_{1}v_{1}u_{6}v_{6}u_{4}v_{4}u_{2}v_{2}u_{7}v_{7}u_{5}v_{5}u_{3}v_{3}u_{1}$, \\
$v_{7}u_{1}v_{3}u_{4}v_{6}u_{7}v_{2}u_{3}v_{5}u_{6}v_{1}u_{2}v_{4}u_{5}v_{7}$ \\
and 
$u_{i}v_{i}u_{i+1}v_{i+1}u_{i+2}v_{i+4}u_{i+5}v_{i+5}u_{i+6}v_{i+6}u_{i+4}v_{i+3}u_{i+3}v_{i+2}u_{i}$, \\
$u_{i+6}v_{i+6}u_{i}v_{i}u_{i+5}v_{i+5}u_{i+3}v_{i+2}u_{i+2}v_{i+4}u_{i+4}v_{i+3}u_{i+1}v_{i+1}u_{i+6}$ and\\
$v_{i+6}u_{i}v_{i+2}u_{i+2}v_{i+1}u_{i+6}v_{i+5}u_{i+3}v_{i+3}u_{i+1}v_{i}u_{i+5}v_{i+4}u_{i+4}v_{i+6}$ \\
for $i=1,2,3,4,5,6,7$. 

Then by Lemma \ref{lemma-HG} (2) we see that counting cycles for a particular edge $e$ of $HG$ will show (1). Similarly we have (2), (3) and (4) by Lemma \ref{lemma-HG} (3), (4) and (5) respectively. 
$\Box$

\vskip 5mm

Summarizing the statements in Lemma \ref{lemma-HG2} we have Table \ref{table-Heawood-cycles}.

\begin{table}[htb]
\centering
\caption{number of $k$-cycles of $HG$}
\scalebox{0.73}
{
\begin{tabular}{|c|c|c|c|c|c|c|c|c|}  \hline
\multirow{2}{*}{$k$} & 
\multirow{2}{*}{$|\Gamma_{k}(HG)|$} & 
\multirow{2}{*}{$|\Gamma_{k}(HG)|\cdot k$} & 
\multirow{2}{*}{$\alpha_{k}(e,HG)$} & 
\multicolumn{3}{c}{$\alpha_{k}(d\cup e,HG)$} & 
\multicolumn{2}{|c|}{$\beta_{k}(d\cup e,HG)$} 
\\ \cline{5-9}
& & & & $d(d,e)=0$ & $d(d,e)=1$ & $d(d,e)=2$ &
$d(d,e)=1$ & $d(d,e)=2$ 

\\ \hline
    $6$ & $28$ & $168$ & $8$ & $4$ & $2$ & $1$ & $2$ & $1$ \\ \hline
    $8$ & $21$ & $168$ & $8$ & $4$ & $2$ & $3$ & $2$ & $1$ \\ \hline
    $10$ & $84$ & $840$ & $40$ & $20$ & $18$ & $17$ & $10$ & $5$ \\ \hline
    $12$ & $56$ & $672$ & $32$ & $16$ & $16$ & $20$ & $8$ & $4$ \\ \hline
    $14$ & $24$ & $336$ & $16$ & $8$ & $12$ & $10$ & $4$ & $2$ \\ \hline
\end{tabular}
}
  \label{table-Heawood-cycles}
\end{table}

\vskip 5mm

\noindent{\bf Proof of Theorem \ref{theorem-Heawood}.} Let $g:HG\to{\mathbb R}^{2}$ be a plane generic immersion illustrated in Figure \ref{HG}. 
We see that $g(HG)$ contains $7$ disjoint crossings between distance $1$ edges and $7$ disjoint crossings between distance $2$ edges. 
Since $\alpha_{6}(d\cup e,HG)=2$ for edges $d$ and $e$ of $HG$ with $d(d,e)=1$ and $\alpha_{6}(d\cup e,HG)=1$ for edges $d$ and $e$ of $HG$ with $d(d,e)=2$, we have
\[
\sum_{\gamma\in\Gamma_{6}(HG)}c(g(\gamma))=2\cdot7+7=21\equiv 1\pmod 2.
\]
Similarly we have
\[
\sum_{\gamma\in\Gamma_{8}(HG)}c(g(\gamma))=2\cdot7+3\cdot7=35\equiv 1\pmod 2
\]
and
\[
\sum_{\gamma\in\Gamma_{10}(HG)}c(g(\gamma))=18\cdot7+17\cdot7=245\equiv 1\pmod 2.
\]
Next we will check that the conditions (1), (2), (3) and (4) of Proposition \ref{proposition-crossing2} for $G=HG$, $\Lambda=\Gamma_{6}(HG)$, $\Gamma_{8}(HG)$ or $\Gamma_{10}(HG)$ and $m=2$ hold. 
Since $\alpha_{6}(e,HG)=\alpha_{8}(e,HG)=8$ and $\alpha_{10}(e,HG)=40$ for every edge $e$ of $HG$, we see that (1) holds. 
For $m=2$, (2) automatically holds. 
Since $\Lambda\subset\Gamma(HG)$ (3) holds for $m=2$ by Remark \ref{remark-crossing}. 
Let $d$ and $e$ be mutually adjacent edges of $HG$. 
Then $\alpha_{6}(d\cup e,HG)=\alpha_{8}(d\cup e,HG)=4$ and $\alpha_{10}(d\cup e,HG)=20$. 
Therefore (4) holds. 
Then by Proposition \ref{proposition-crossing2} we have 
\[
\sum_{\gamma\in\Gamma_{6}(HG)}c(f(\gamma))\equiv \sum_{\gamma\in\Gamma_{6}(HG)}c(g(\gamma))\pmod 2,
\]
\[
\sum_{\gamma\in\Gamma_{8}(HG)}c(f(\gamma))\equiv \sum_{\gamma\in\Gamma_{8}(HG)}c(g(\gamma))\pmod 2
\]
and
\[
\sum_{\gamma\in\Gamma_{10}(HG)}c(f(\gamma))\equiv \sum_{\gamma\in\Gamma_{10}(HG)}c(g(\gamma))\pmod 2.
\]
Therefore we have 
\[
\sum_{\gamma\in\Gamma_{6}(HG)}c(f(\gamma))\equiv\sum_{\gamma\in\Gamma_{8}(HG)}c(f(\gamma))\equiv\sum_{\gamma\in\Gamma_{10}(HG)}c(f(\gamma))\equiv 1\pmod 2. 
\]
We note that $\alpha_{12}(e,HG)=32$ for every edge $e$ of $HG$, $\alpha_{12}(d\cup e,PG)=16$ for any adjacent edges $d$ and $e$, $\alpha_{12}(d\cup e,PG)=16$ for any edges $d$ and $e$ with $d(d,e)=1$ and $\alpha_{12}(d\cup e,PG)=20$ for any edges $d$ and $e$ with $d(d,e)=2$. 
Then by Proposition \ref{proposition-crossing1} we have 
\[
\sum_{\gamma\in\Gamma_{12}(HG)}c(f(\gamma))\equiv 0\pmod 4.
\]
We note that $\alpha_{14}(e,HG)=16$ for every edge $e$ of $HG$, $\alpha_{14}(d\cup e,PG)=8$ for any adjacent edges $d$ and $e$, $\alpha_{14}(d\cup e,PG)=12$ for any edges $d$ and $e$ with $d(d,e)=1$ and $\alpha_{14}(d\cup e,PG)=10$ for any edges $d$ and $e$ with $d(d,e)=2$. 
Then by Proposition \ref{proposition-crossing1} we have 
\[
\sum_{\gamma\in\Gamma_{14}(HG)}c(f(\gamma))\equiv 0\pmod 2.
\]
$\Box$

\section{Rotation numbers of cycles in a plane immersed graph}\label{rotation-numbers} 

In this section we give a proof of Theorem \ref{theorem-zero-rotation}.

\begin{Lemma}\label{lemma-zero-rotation}

Let $G$ be a $2$-connected loopless graph. Let $u$ and $v$ be distinct vertices of $G$ such that the graph obtained from $G$ by adding an edge joining $u$ and $v$ does not have $K_{4}$ as a minor. 
Then there exists a plane generic immersion $f:G\to{\mathbb R}^{2}$ with the following properties. 

\begin{enumerate}

\item[{\rm (1)}]

It holds that ${\rm rot}(f(\gamma))=0$ for every cycle $\gamma$ of $G$. 

\item[{\rm (2)}]

Let $h:{\mathbb R}^{2}\to{\mathbb R}$ be a height function defined by $h(x,y)=y$. 
Then $h(f(G))=[h(f(v)),h(f(u))]$. If $P$ is a path of $G$ joining $u$ and $v$, then $h\circ f$ maps $P$ homeomorphically onto the closed interval $[h(f(v)),h(f(u))]$. 

\end{enumerate}

\end{Lemma}

\noindent{\bf Proof.} The following is a proof by induction on the number of vertices of $G$. 

\begin{figure}[htbp]
      \begin{center}
\scalebox{0.6}{\includegraphics*{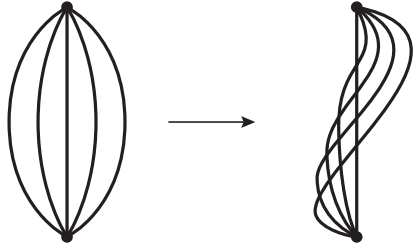}}
      \end{center}
   \caption{All cycles have rotation number $0$}
  \label{theta-n}
\end{figure}

First suppose $|V(G)|=2$. Then $G$ is a $\theta_{n}$-curve graph and a plane generic immersion described in Figure \ref{theta-n}, where the case $n=5$ is illustrated, satisfies (1) and (2). 

Let $k$ be a natural number greater than or equal to $2$. Suppose that the claim is true if $|V(G)|\leq k$. Suppose $|V(G)|=k+1$. 
Let $X_{1},\cdots,X_{n}$ be the connected components of the topological space $G\setminus\{u,v\}$ and $H_{1},\cdots,H_{n}$ the closures of them in $G$. Namely $H_{i}=X_{i}\cup\{u,v\}$ and $H_{i}$ is regarded as a subgraph of $G$ for $i=1,\cdots,n$. 
Let $H_{i}'$ be a graph obtained from $H_{i}$ by adding an edge $e_{i}$ joining $u$ and $v$. By the assumption on $G$ we see that $H_{i}'$ is a $2$-connected loopless graph that does not have $K_{4}$ as a minor. 

\begin{figure}[htbp]
      \begin{center}
\scalebox{0.6}{\includegraphics*{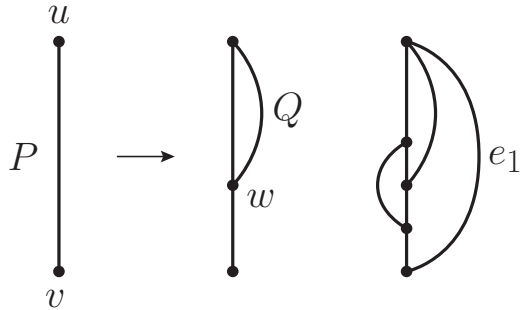}}
      \end{center}
   \caption{$n$ is greater than $1$}
  \label{path}
\end{figure}

We will show that $n$ is greater than $1$. Suppose to the contrary that $n=1$. Then $G=H_{1}$. Let $P$ be a path of $H_{1}$ joining $u$ and $v$. Suppose that there exists another path $Q$ of $H_{1}$ joining $u$ and $v$ such that $P\cap Q=\{u,v\}$. Since $X_{1}=H_{1}\setminus\{u,v\}$ is connected there is a path $R$ of $H_{1}$ away from $u$ and $v$ joining a vertex of $P$ and a vertex of $Q$. Then we see that $H_{1}'$ has $K_{4}$ as a minor. 
Therefore there are no such path of $H_{1}$. 
Since $H_{1}=G$ is $2$-connected, $H_{1}$ has no cut vertices. Therefore ${\rm deg}(u)\geq2$. 
Then there exists a path $Q$ of $H_{1}$ joining $u$ and a vertex $w$ of $P$ such that $P\cap Q=\{u,w\}$. 
We may suppose that $w$ is closest to $v$ on $P$ among all choices of such path $Q$. Then we have either $w$ is a cut vertex of $H_{1}$ or $H_{1}'$ has $K_{4}$ as a minor. See Figure \ref{path} for the latter case. 
We note that the second vertex from bottom in the right graph of Figure \ref{path} may be $v$. 
Both contradict to the assumption. Thus we have $n\geq2$. 

Since $|V(G)|\geq3$ at least one of $H_{1},\cdots,H_{n}$, say $H_{1}$, contains at least $3$ vertices. 

\begin{figure}[htbp]
      \begin{center}
\scalebox{0.6}{\includegraphics*{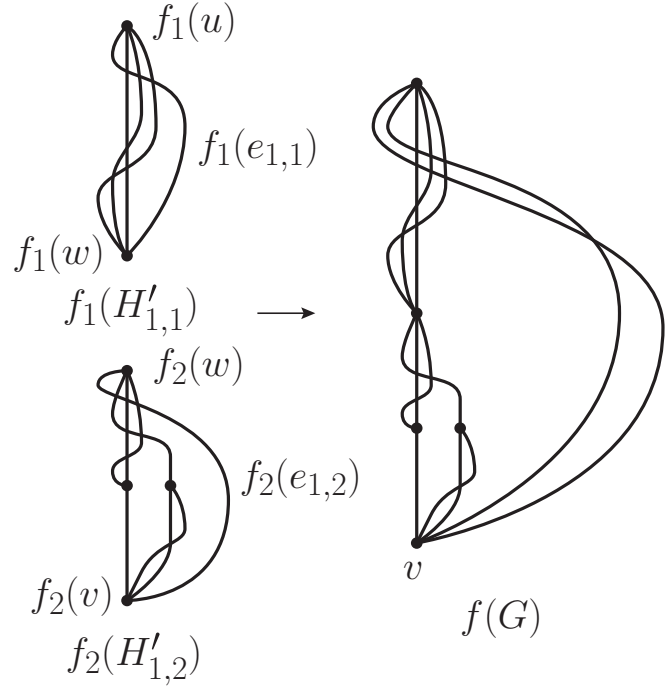}}
      \end{center}
   \caption{A construction}
  \label{zero-immersion}
\end{figure}

Suppose that all other $H_{i}$ contains exactly two vertices. Then $H_{i}$ contains exactly one edge that joins $u$ and $v$ for $i=2,\cdots,n$. 
Then by an argument similar to the previous one we see that $H_{1}$ has a cut vertex $w$. Let $H_{1,1}$ and $H_{1,2}$ be connected subgraphs of $H_{1}$ such that $H_{1}=H_{1,1}\cup H_{1,2}$ and $H_{1,1}\cap H_{1,2}=\{w\}$. We may suppose without loss of generality that $H_{1,1}$ contains $u$ and $H_{1,2}$ contains $v$. Let $H_{1,1}'$ be a graph obtained from $H_{1,1}$ by adding an edge $e_{1,1}$ joining $u$ and $w$. Let $H_{1,2}'$ be a graph obtained from $H_{1,2}$ by adding an edge $e_{1,2}$ joining $w$ and $v$. 
Then $H_{1,1}'$ and $u$ and $w$ satisfy the induction hypothesis and $H_{1,2}'$ and $w$ and $v$ also satisfy the induction hypothesis. 
Let $f_{1}:H_{1,1}'\to{\mathbb R}^{2}$ and $f_{2}:H_{1,2}'\to{\mathbb R}^{2}$ be plane generic immersions that satisfy (1) and (2). 
By a translation we may assume $f_{1}(w)=f_{2}(w)$. 
Let $f_{0}:H_{1}\to{\mathbb R}^{2}$ be a plane generic immersion defined by $f_{0}|_{H_{1,1}}=f_{1}|_{H_{1,1}}$ and $f_{0}|_{H_{1,2}}=f_{2}|_{H_{1,2}}$. Then by a construction similar to that illustrated in Figure \ref{theta-n} we have a plane generic immersion $f:G\to{\mathbb R}^{2}$ with $f|_{H_{1}}=f_{0}$ that satisfies (1) and (2). 
See for example Figure \ref{zero-immersion} where the case $n=3$ is illustrated. 

\begin{figure}[htbp]
      \begin{center}
\scalebox{0.6}{\includegraphics*{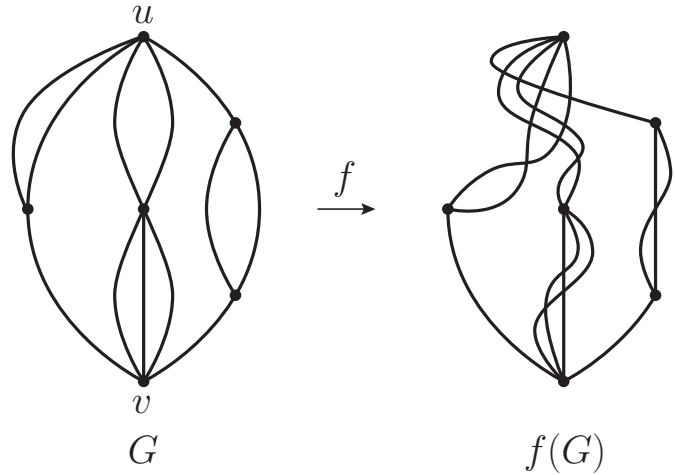}}
      \end{center}
   \caption{An example}
  \label{zero-immersion2}
\end{figure}

Next suppose that another $H_{i}$, say $H_{2}$, contains at least three vertices. Then by a construction similar to that illustrated in Figure \ref{zero-immersion}, we have a plane generic immersion of $G$ that satisfies (1) and (2). 
See for example Figure \ref{zero-immersion2}. This completes the proof.
$\Box$


\vskip 5mm

\noindent{\bf Proof of Theorem \ref{theorem-zero-rotation}.} Suppose that $K_{4}$ is a minor of $G$. Since $K_{4}$ is $3$-regular, $G$ contains a subgraph that is homeomorphic to $K_{4}$. 
Then by Corollary \ref{corollary-K4} we see that every plane generic immersion of $G$ contains a cycle with non-zero rotation number. Therefore (1) implies (2). 
Suppose that $G$ does not have $K_{4}$ as a minor. We will show that there exists a plane generic immersion $f:G\to{\mathbb R}^{2}$ such that ${\rm rot}(f(\gamma))=0$ for every cycle $\gamma$ of $G$. 
By considering the block decomposition, it is sufficient to show the case that $G$ is $2$-connected and loopless. Then by Lemma \ref{lemma-zero-rotation} we have such a plane generic immersion. This completes the proof. 
$\Box$

\section{The total Thurston-Bennequin number of a Legendrian embedding of a finite graph}\label{total-Thurston-Bennequin-number} 

\begin{Lemma}\label{lemma-total-tb} Let $G$ be a finite graph and $f:G\to {\mathbb R}^{3}$ a Legendrian embedding of $G$. 
Let $j$ and $k$ be natural numbers. Suppose that there exists a rational number $q$ such that the following holds. 
\begin{enumerate}

\item[{\rm (1)}]

For any edge $e$ of $G$, $\alpha_{k}(e,G)=q\cdot\alpha_{j}(e,G)$, 

\item[{\rm (2)}]

For any pair of mutually adjacent edges $d$ and $e$ of $G$, $\alpha_{k}(d\cup e,G)=q\cdot\alpha_{j}(d\cup e,G)$, 

\item[{\rm (3)}]

For any pair of mutually disjoint oriented edges $d$ and $e$ of $G$, $\beta_{k}(d\cup e,G)=q\cdot\beta_{j}(d\cup e,G)$. 

\end{enumerate}
Then
\[
TB_{k}(f)=q\cdot TB_{j}(f).
\]
\end{Lemma}

\noindent{\bf Proof.} Let $\gamma$ be a cycle of $G$. 
Let $D=D(f(\gamma))$ be a Lagrangian projection of a Legendrian knot $f(\gamma)$. 
Let $w(D)$ be the writhe of $D$. 
It is known that $\displaystyle{tb(f(\gamma))=w(D)}$. 
We note that each crossing of $D$ is a self crossing, an adjacent crossing or a disjoint crossing. 
A self crossing of an edge $e$ of $G$ contribute $\pm\alpha_{j}(e,G)$ to $TB_{j}(f)$ and $\pm\alpha_{k}(e,G)$ to $TB_{k}(f)$. 
Similarly an adjacent crossing of edges $d$ and $e$ contribute $\pm\alpha_{j}(d\cup e,G)$ to $TB_{j}(f)$ and $\pm\alpha_{k}(d\cup e,G)$ to $TB_{k}(f)$. A disjoint crossing of edges $d$ and $e$ contribute $\pm\beta_{j}(d\cup e,G)$ to $TB_{j}(f)$ and $\pm\beta_{k}(d\cup e,G)$ to $TB_{k}(f)$. 
By taking the sum of these numbers we have the result. 
$\Box$

\vskip 5mm

\noindent{\bf Proof of Theorem \ref{theorem-Petersen-tb}.} By Lemma \ref{lemma-total-tb} and Lemma \ref{lemma-PG2} we have the result. 
$\Box$

\vskip 5mm

\noindent{\bf Proof of Theorem \ref{theorem-Heawood-tb}.} By Lemma \ref{lemma-total-tb} and Lemma \ref{lemma-HG2} we have the result. 
$\Box$

\vskip 3mm

\section*{Acknowledgments} The authors are grateful to the referee for his/her helpful comments.

\vskip 3mm


{\normalsize

}

\end{document}